\documentclass[10pt]{article} 
 \usepackage[margin=1in]{geometry}
\usepackage[sectionbib]{natbib}
\usepackage{array,epsfig,fancyheadings,rotating}
\usepackage{sectsty}
\usepackage{setspace}
\usepackage{titlesec}
\usepackage{natbib}
\usepackage[]{hyperref}  
\usepackage[dvipsnames]{xcolor}
\usepackage{etoolbox}
\usepackage[font=small,skip=-2pt]{caption}
\usepackage{amsmath}
\usepackage{amssymb}
\usepackage{amsfonts}
\usepackage{multirow}
\usepackage{amsthm}
\usepackage{adjustbox}
\usepackage{thm-restate}
\usepackage{enumitem}
\usepackage{bbm}
\usepackage{algorithmic}
\usepackage{amsmath,amssymb,bm,mathrsfs,mathtools}
\usepackage{etoolbox}
\usepackage[plain,noend,linesnumbered]{algorithm2e}
\AtBeginEnvironment{quote}{}

\newtheorem{theorem}{Theorem}[section]

\newtheorem{assumption}[theorem]{Assumption}
\newtheorem{proposition}[theorem]{Proposition}
\newtheorem{lemma}[theorem]{Lemma}
\newtheorem{corollary}[theorem]{Corollary}
\newtheorem{remark}[theorem]{Remark}
\theoremstyle{definition}

\newtheorem{example}[theorem]{Example}

\def\hC{{\widehat C}}

\def\bE {{\bf E}}

\def\mF{{\mathcal F}}
\def\mP{{\mathcal P}}
\def\mV{{\mathcal V}}
\def\mZ{{\mathcal Z}}
\def\mX{{\mathcal X}}
\def\bP{{\mathbbm{P}}}
\def\bE{{\mathbbm{E}}}

\def\eps {{\varepsilon}}
\newcommand{\indep}{\perp \!\!\! \perp}

\newcommand{\real}{\mathbb{R}}

\numberwithin{equation}{section}

\allowdisplaybreaks
\title{Localized Conformal Prediction: \\
A Generalized Inference Framework to Conformal Prediction}
\author{Leying  Guan\thanks{Dept. of  Biostatistics,
    Yale University, leying.guan@yale.edu}
    }
 \date{}
\begin{document}
\maketitle

\begin{abstract}
We propose a new inference framework called localized conformal prediction. It generalizes the framework of conformal prediction by offering  a single-test-sample adaptive construction that emphasizes a local region around this test sample, and can be combined with different conformal score constructions. The proposed framework enjoys an assumption-free finite sample marginal coverage guarantee, and it also offers additional local coverage guarantees under suitable assumptions. We demonstrate how to change from conformal prediction to localized conformal prediction using several conformal scores, and we illustrate a potential gain via numerical examples.
\end{abstract}

\section{Introduction}
\label{sec:intro}
Conformal prediction (CP) is an increasingly popular framework for measuring prediction uncertainty. Let $Z_i\coloneqq (X_i, Y_i)$ for $i=1,\ldots,n$ be i.i.d.\ regression data from some joint distribution $\mP_{XY}$, where $X_i\in \real^p$ is the feature and $Y_i\in\real$ is the response. Given a new feature $X_{n+1}$ with its response $Y_{n+1}$ unobserved, the goal of CP is to construct a prediction interval (PI) $ C(X_{n+1})$ that covers  $Y_{n+1}$ with probability at least $\alpha$:
\begin{equation}
\label{eq:goal}
\bP(Y_{n+1}\in C(X_{n+1})) \geq \alpha,
\end{equation}
for some desired coverage level $\alpha\in (0,1)$ (usually close to 1).  Setting $Z_{n+1}=(X_{n+1}, Y_{n+1})$ as the $(n+1)^\text{th}$ observation, CP achieves (\ref{eq:goal}) under the assumption that $Z_{n+1}$ is also independently generated from $\mP_{XY}$, without additional distributional assumptions on $\mP_{XY}$ itself \citep{vovk2005algorithmic, shafer2008tutorial, vovk2009line,lei2014distribution, lei2018distribution}.  

Let $\mZ = \{Z_1,\ldots,Z_{n+1}\}$ be the unordered set of feature-response pairs, including $Z_{n+1}=(X_{n+1},Y_{n+1})$. CP is based on a conformal score function $V(.)$ for observations $z=(x,y)$, whose form may depend also on the unordered data $\mZ$, i.e. $V(z) = V(z;\mZ)$. We will consider score functions $V(.)$ where large values of $V(z)$ indicate that $z$ is less likely to be a sample from $\mP_{XY}$. For instance, we may choose $V(z)=|\hat\mu(x) - y|$ for a prediction function $\hat\mu(.)$ that is learned from the data $\mZ$ or from a separate independent data set.

It is guaranteed that $V_i \coloneqq V(Z_i)$ are exchangeable when $Z_1,\ldots,Z_{n+1}$ are i.i.d. Thus, letting $Q(\alpha; V_{1:{n+1}})$ denote the level-$\alpha$ quantile of the empirical distribution of $V_1,\ldots,V_{n+1}$,  we have
\begin{equation}
\label{eq:eq1}
\bP\left\{V_{n+1} \leq Q(\alpha; V_{1:n+1})\right\} \geq \alpha.
\end{equation}
CP constructs the level-$\alpha$ PI $C(X_{n+1})$ for $Y_{n+1}$ by inverting the above relationship for $V_{n+1}$:
\begin{equation}
\label{eq:eq2}
 C(X_{n+1}) = \{y: V(X_{n+1},y)\leq  Q(\alpha; V_{1:n}\cup V(X_{n+1},y))\}.
\end{equation}
Note that if the form of $V(z)=V(z;\mZ)$ depends also on $\mZ$, then in the above, each score $V_i=V(Z_i;\mZ)$ depends also on $y$ and is understood to be evaluated at $Z_{n+1}=(X_{n+1},y)$. By the guarantee (\ref{eq:eq1}), $C(X_{n+1})$ constructed in this way satisfies (\ref{eq:goal}) for any distribution $\mP_{XY}$.

It is common in data applications for the conditional distribution of $Y$ given $X=x$ to be heterogeneous across different values of $x \in \real^p$. In such settings, it is desirable for the constructed PI $C(X_{n+1})$ to adapt to this heterogeneity.  However, by definition, the CP interval $C(X_{n+1})$ is based on the global exchangeability of the conformal scores $V_1,\ldots,V_{n+1}$, and depends equally on scores where $X_i$ is far from $X_{n+1}$ as on scores where $X_i$ is close to $X_{n+1}$. To adapt to the heterogeneity of $Y$ given $X=x$,  one active area of research has been to design the score function $V(.)$ to directly capture this heterogeneity, in a way so that the quantiles of $V(.)$ are more homogeneous across different $x \in \real^p$ \citep{lei2014distribution, izbicki2019flexible, lei2018distribution, romano2019conformalized, gupta2021nested}. For example, in \cite{romano2019conformalized}, the authors consider the quantile regression score $V(z) =\max\{\hat q_{lo}(x)-y,y-\hat q_{hi}(x)\}$ where $\hat q_{lo}(x)$ and $\hat q_{hi}(x)$ are estimated quantiles for the conditional distribution of $Y$ given $X=x$. However, this approach may yield deteriorated performance when these quantile functions are difficult to estimate for some regions of the feature space.

In this paper, we take a different approach, and generalize the inference framework itself by weighting the conformal scores $V_1,\ldots,V_n$ differently based on the observed feature value $X_{n+1}$. Our method places more weight on scores $V_i$ for which $X_i$ belongs to a local region around $X_{n+1}$. Performing conformal inference while emphasizing the unique role of $X_{n+1}$ is an interesting and open problem, and we provide the first such generalization with theoretical guarantees. We call this generalized framework localized conformal prediction (LCP), which can be flexibly combined with recently developed conformal score functions. 

The main idea of LCP is to introduce a localizer around $X_{n+1}$, and up-weight samples close to $X_{n+1}$ according to this localizer. For example, we may take the localizer $H(X_{n+1},X_i) = e^{-5|X_i - X_{n+1}|}$, consider the weighted empirical distribution where $V_i$ has weight proportional to $H(X_{n+1},X_i)$, and include the value $y$ in $C(X_{n+1})$ if and only if $V(X_{n+1}, y)$ is smaller than the $\tilde\alpha$ quantile of this weighted distribution. As this weighted distribution is no longer exchangeable, we will need to choose $\tilde\alpha$ strategically  to guarantee finite-sample coverage as described in (\ref{eq:goal}).

We demonstrate the difference between LCP and CP with a simple example: Features $X\sim \operatorname{Unif}(-5,5)$ follow a uniform distribution on $[-5,5]$, and the response $Y$ given $X$ follows a mean-zero normal distribution with heterogeneous variance across $X$:
\[
Y|X \sim\left\{\begin{array}{ll}\cos(\frac{\pi}{10}X_i)\times N(0,1),&\mbox{if } |X|\leq 4.5,\\
 2\times N(0,1),&\mbox{if } |X|> 4.5.
\end{array}\right.
\]
We fix the desired coverage level $\alpha = 0.95$, take $n=1000$ samples, and perform both CP and LCP using two score functions: (1) the regression score $V(z) = |\mu(x)-y|=|y|$ where here $\mu(x)=0$ \citep{lei2018distribution}, and (2) the quantile regression score $V(z) =\max\{\hat q_{lo}(x)-y,y-\hat q_{hi}(x)\}$ \citep{romano2019conformalized} where $\hat q_{lo}$ and $\hat q_{hi}$  are $0.025$ and $0.975$ quantile curves estimated from 2000 independent samples using a neural network model as described in Section \ref{sec:empirical}. The localizer for LCP is $H(X_{n+1},X_i)=e^{-5|X_i-X_{n+1}|}$ as above. We refer to the two corresponding CP procedures as CR and CQR, and the two LCP procedures as LCR and LCQR.

\begin{figure}
\begin{center}
\includegraphics[width = 1\textwidth]{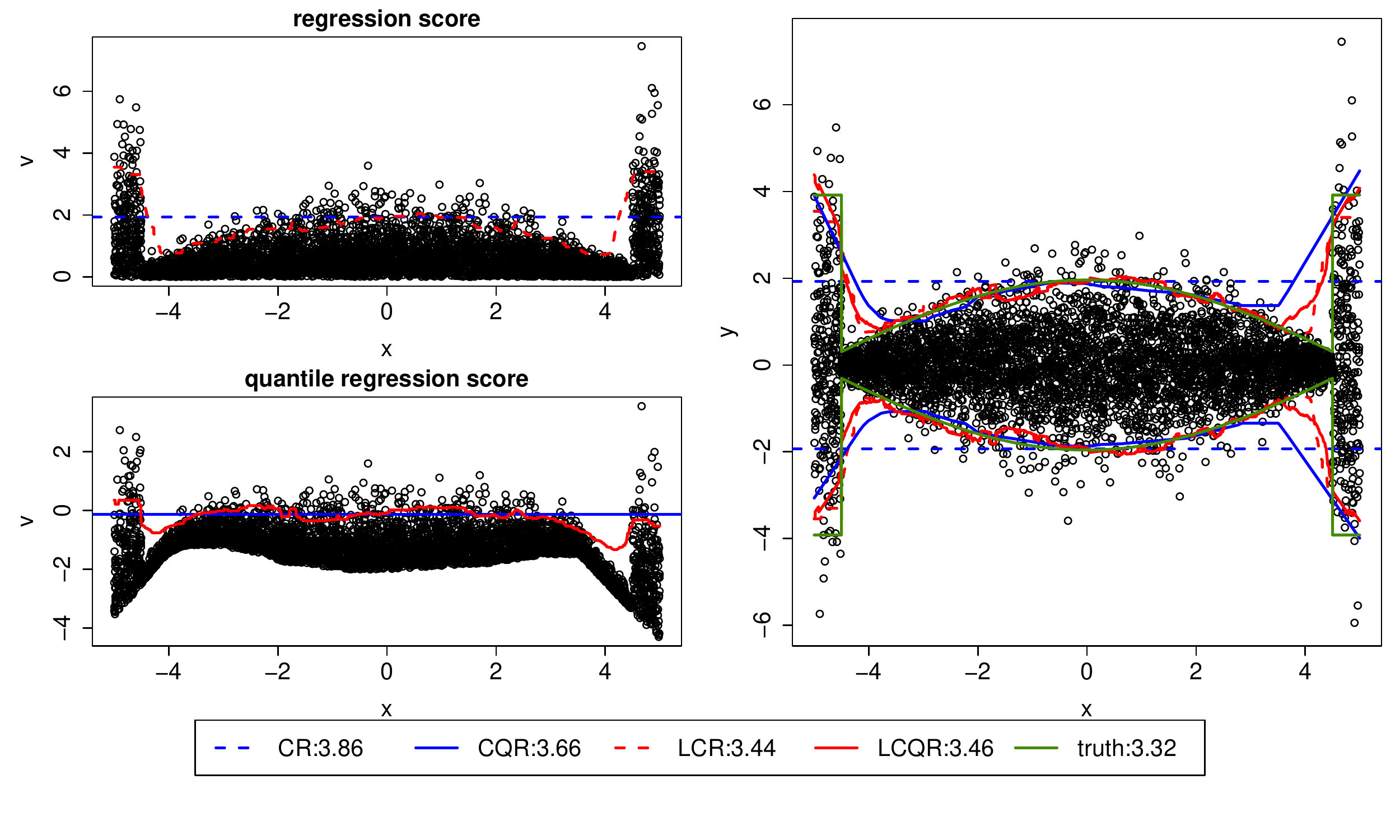} 
\caption{\em Conformal bands (blue) and localized conformal bands (red) using regression score (dashed) and quantile regression score (solid). Bands are shown for the conformal score $V_{n+1}$ (left) and response value $Y_{n+1}$ (right). True prediction bands for the distribution of $Y$ given $X$ are shown on the right in green, and dots show the realized test observation values. The legend indicates average PI length, which is shorter for the localized procedures.}
\label{fig:illustration1}
\end{center}
\vskip -.3 in
\end{figure}

The left panel of Figure \ref{fig:illustration1} shows the conformal confidence bands for $V_{n+1}$ using CR\slash LCR (upper left, blue\slash red, dashed curve) and CQR\slash LCQR (lower left, blue\slash red, solid curve). The right panel shows the inverted PI for $Y_{n+1}$ using the four procedures. The green curves on the right panel represent the true level-$\alpha$ confidence bands for $Y$ given $X$.  This example demonstrates that, by definition, the CR and CQR intervals are homogeneous for $V$. In this example, the CR intervals are furthermore homogeneous for $Y$. CQR provides a heterogeneous PI for $Y$ by inverting the interval for $V$. However,  the true quantile functions are hard to estimate at the two ends here, and thus some heterogeneity of $V_{n+1}$ still remains even for the quantile regression score. In comparison, LCP introduces more flexibility by directly constructing intervals that are heterogeneous for $V$. It yields an improvement even when applied to the quantile regression score, where it better captures the remaining heterogeneity of this score.

We summarize our contributions as follows:
 
 \begin{itemize}
 \item We generalize the probabilistic framework of CP to LCP,  where we assign a unique role to the test point by introducing a localizer around it. The generalized framework still enjoys a distribution-free and finite-sample marginal coverage guarantee. CP is a special case of LCP where the localizer takes a constant value.
 \item We focus on sample-splitting LCP and develop an efficient implementation. We also demonstrate how to combine LCP with some recently developed conformal scores, with numerical examples.
 \item We investigate the local behavior of sample-splitting LCP and show that it enjoys additional local coverage guarantees under proper assumptions.
 \end{itemize}
We postpone  all proofs to Appendix \ref{app:proofmain}  in the online Supplement.

\section{LCP: A generalization of conformal prediction}
\label{sec:LCPexchange}
\subsection{Notations}
\label{subsec:notation1}
For any distribution $\mathcal{F}$ on $\real$, we define its level $\alpha$ quantile as 
\[
Q(\alpha;\mathcal{F}) = \inf\{t: \bP_{T\sim \mF}\{T \leq t\} \geq \alpha\}.
\]
Let $\mX=\{X_1,\ldots, X_{n+1}\}$ be the unordered set of feature values from all $(n+1)$ samples. LCP weights samples differently based on a bi-variate localizer function $H(x, x'):\real^p \times \real^p\mapsto [0,1]$, whose function form may depend on the data through only $\mX$.  We  require $H(x, x) = 1$ for all $x$ and use $H(x, x')$ to capture the dissimilarity between two given feature values. In the Introduction, we considered $H(x, x') = \exp(-5|x-x'|)$ as an example where $H(X_{n+1},X_i)$ is the localizer evaluated at $X_{n+1}$ and $X_i$. The localizer function $H(x,x')$ is  used on construct different weighted distributions for performing LCP. Define $H_i(.)\coloneqq H(X_i, .)$ as the localizer centered at $X_i$, and $H_{i,j} \coloneqq H_i(X_j) = H(X_i, X_j)$ as a measure of dissimilarity between samples $X_i$ and $X_j$.    

Let $\delta_v$ be a point mass at $v\in \real$. Define weighted distributions 
\begin{align*}
&\hat \mF_i \coloneqq \left(\sum^{n+1}_{j=1} p^H_{i,j}\delta_{V_j}\right),\; \mbox{for all } i = 1,\ldots, n+1,
\end{align*}
where the empirical weights $p^{H}_{i,j} \coloneqq \frac{H_{ij}}{\sum^{n+1}_{k=1} H_{ik}}$ for $j = 1,\ldots, n+1$ are constructed using  the localizer  centered at $X_i$.  We also define
\begin{align*}
&\hat \mF \coloneqq \left(\sum^{n}_{j=1} p^H_{n+1,j}\delta_{V_j}+p^H_{n+1,n+1}\delta_\infty\right),
\end{align*}
as the distribution when replacing $V_{n+1}$ by $\infty$ in $\hat\mF_{n+1}$.  Note that both $V_{n+1}$ and $V_i$ for $i=1,\ldots, n$ may depend on $Y_{n+1}$ when $V(.)$ depends on the set $\mZ$. Consequently, there could be a dependence on $Y_{n+1}$ from both $\hat\mF$ and $\hat\mF_i$ for $i=1,\ldots, n+1$. We have masked such a dependence for  convenience. 

Throughout this paper, we call $\{Z_1,\ldots, Z_n\}$ the calibration set and assume $Z_i\overset{i.i.d}{\sim} \mP_{XY}$ for $i=1,\ldots, n+1,$ and $\alpha\in (0,1)$ is a constant and user-specified targeted coverage. 

\subsection{LCP and marginal coverage guarantee}
We now establish the probabilistic guarantees of LCP regarding its marginal coverage. Instead of using the level $\alpha$ quantile of the empirical distribution as in CP,  LCP considers a level $\tilde\alpha$ quantile of a weighted empirical distribution, with weight proportional to $H_{n+1,i}$. Recall that $H_{n+1,i}$ measures the distance between a training sample $X_i$ and the test sample $X_{n+1}$. This weighted distribution allows more emphasis on training samples closer to $X_{n+1}$.  

Theorem \ref{thm:general1} states how we can choose $\tilde\alpha$ to achieve finite sample coverage.  In Theorem \ref{thm:general2}, we show that a randomized decision rule can lead to a PI with exact coverage.  Let $\Gamma =  \{ \sum_{k\in I_i}p^{H}_{ik}: i=1,\ldots,n+1, I_i \subseteq \{1,\ldots, n+1\}\}$ to represent all possible empirical CDF function values from weighted distributions $\hat\mF_i$ for $i=1,\ldots, n+1$, under all possible ordering of $V_1,\ldots, V_{n+1}$.
\begin{theorem}
\label{thm:general1} 
Let $\tilde\alpha$ be the  smallest value in $\Gamma$ such that
\begin{equation}
\label{eq:goal1general}
\sum^{n+1}_{i=1}\frac{1}{n+1}\mathbbm{1}_{V_i \leq Q(\tilde\alpha; \hat\mF_i)}  \geq \alpha  
\end{equation}
Then $\bP\left\{V_{n+1} \leq Q(\tilde\alpha;\hat \mF_{n+1} ) \right\}\geq \alpha$. Equivalently, $\bP\left\{V_{n+1} \leq Q(\tilde\alpha;\hat \mF) \right\}\geq \alpha$.
\end{theorem} 
\begin{remark}
\label{remark1}
If  $H_{i,j} = 1$ for all $i, j = 1,\ldots, n+1$, then we have $\hat \mF=\frac{1}{n+1}\left( \sum_{i=1}^n\delta_{V_i}+\delta_{\infty}\right)$ and $\hat\mF_i =\frac{1}{n+1}\left( \sum_{i=1}^{n+1}\delta_{V_i}\right)$ for all $i=1,\ldots, n+1$. Then (\ref{eq:goal1general}) holds  if and only if $\tilde\alpha\geq \alpha$ by definition. Also, $\Gamma = \{\frac{k}{n+1}:k=1,\ldots, n+1\}$. Thus,  we recover usual conformal prediction \citep{vovk2005algorithmic}, and
\[
\bP\left\{V_{n+1} \leq Q\left(\frac{k}{n+1};  \frac{1}{n+1}\sum_{i=1}^{n+1}\delta_{V_i}\right) \right\}\geq \alpha, \; \mbox{for }k \geq \lceil (n+1)\alpha\rceil.
\]
\end{remark}

Here, we provide some intuition for why such $\tilde\alpha$ can guarantee level $\alpha$ coverage. Conformal prediction relies on the exchangeability of data. Conditional on the set  $\mZ$,  the set of observed values $ \mV=\{v_1, \ldots, v_{n+1}\}$ for $V_{1:(n+1)}$ is fixed and $V_{n+1}$ has equal probability of taking each of value in  $\mV$. Hence $Q(\alpha, \frac{1}{n+1}\sum_{i=1}^{n+1}\delta_{v_{i}})$ leads to a  coverage guarantee conditional on the observed values, and a marginal coverage guarantee after marginalizing over all value sets. When our PI is constructed as $Q(\tilde\alpha, \hat\mF_{n+1})$, since $\hat\mF_{n+1}$ changes as we permute the value assignments, we need to account for this change  while calculating the conditional coverage. The left-hand-side of (\ref{eq:goal1general}) turns out to be this coverage conditional on $\{v_1, \ldots, v_{n+1}\}$ for any given $\tilde\alpha$.  As in CP, we can invert the relationship (\ref{eq:goal1general}) to construct the PI for $Y_{n+1}$.
\begin{corollary}
\label{cor:generalC} 
In the setting of Theorem \ref{thm:general1}, define $\tilde\alpha(y)\in \Gamma$ as the smallest value in $\Gamma$ such that (\ref{eq:goal1general}) holds at $Z_{n+1}=(X_{n+1},y)$. Let $C(X_{n+1}) \coloneqq \{y:V_{n+1} \leq Q(\tilde\alpha(y);\hat \mF)\}$. Then, we have $\bP\left\{Y_{n+1}\in  C(X_{n+1})\right\}\geq \alpha$.
\end{corollary}

What will happen if we simply let $\tilde\alpha = \alpha$ without tuning it based on (\ref{eq:goal1general})? The answer depends on the localizer $H$.    Setting $\tilde\alpha=\alpha$ can lead to over-coverage in the simple example described by Proposition  \ref{exm:counter1} below, where we tend to assign too little weight to the calibration samples. A more interesting example is given in Proposition \ref{exm:counter2} below, showing that we may end up achieving arbitrarily bad under-coverage by naively setting $\tilde\alpha=\alpha$.
\begin{proposition}
\label{exm:counter1}
Consider the localizer $H(x_1, x_2) =\exp(-\frac{|x_1-x_2|}{h})$ with some small $h>0$, such that $\bP(\sum^{n+1}_{i=1} H(X_{n+1}, X_i) < \frac{1}{1-\alpha})\geq \varepsilon\in (\alpha,1)$. Then,  
\[
\bP(Q(\alpha;\hat\mF)=\infty)\geq \varepsilon, \qquad \bP(V_{n+1}\in C(X_{n+1}))\geq \varepsilon.
\]
\end{proposition}
\begin{proposition}
\label{exm:counter2}
Let $\{e_j:j = 1,\ldots, p\}$ be the standard basis in $\real^p$. Set $q_1 =\frac{(1-\alpha)}{2p(1-\alpha)+\alpha}$ and $q_0 = \frac{\alpha}{2p(1-\alpha)+\alpha}$. Suppose that the feature $X\in \real^p$ and response $Y$ are  distributed as 
\[
Y|X \sim   \left\{\begin{array}{ll}  \operatorname{Unif}(-1, 1),&\mbox{when } X\neq 0.\\
0,&\mbox{otherwise} .
\end{array}\right., \;X = \left\{\begin{array}{ll}e_j,&\;\mbox{w.p. } q_1 ,\; \mbox{for all } j =1,\ldots, p,\\
-e_j,&\;\mbox{w.p. }q_1,\; \mbox{for all } j =1,\ldots, p,\\
0,&\;\mbox{w.p. }q_0.
\end{array}\right.
\]
Let $V(Z_i)=|Y_i|$ be the regression score. Then for any constant $p \geq 1$, we have $\lim_{n\rightarrow\infty}\bP(V_{n+1}\leq Q(\alpha; \hat\mF)) =q_0<\alpha$.
\end{proposition}

Proposition \ref{exm:counter2} shows an example in which we no longer enjoy the distribution-free marginal coverage guarantee, and the under-coverage can be arbitrarily poor for large $p$. Hence, strategically choosing $\tilde\alpha$ is crucial to obtain such a guarantee. We note that this distribution-free marginal coverage guarantee is usually motivation for using conformal prediction as opposed to other model-based prediction intervals.

As in the case of CP,  we may not have exact level $\alpha$-coverage due to rounding issues using a non-random construction rule. However, we can have exact $\alpha$-coverage if we allow for some additional randomness,  as stated in Theorem \ref{thm:general2}.
\begin{theorem}
\label{thm:general2} 
Consider  the setting of  Theorem \ref{thm:general1}. Let $\tilde\alpha_1\slash\tilde\alpha_2$ be the smallest\slash largest value in $\Gamma\cup\{0\}$ such that 
\[
\alpha_1\coloneqq \sum^{n+1}_{i=1}\frac{1}{n+1}\mathbbm{1}_{V_i \leq Q(\tilde\alpha_1;\hat\mF_i)} \geq \alpha,\qquad \alpha_2\coloneqq \sum^{n+1}_{i=1}\frac{1}{n+1}\mathbbm{1}_{V_i \leq Q(\tilde\alpha_2;\hat\mF_i)} < \alpha.
\]
Set $\tilde\alpha = \left\{\begin{array}{cc}\tilde \alpha_1&w.p.\;\frac{\alpha - \alpha_2}{\alpha_1-\alpha_2}\\
\tilde\alpha_2 & w.p.\;\frac{\alpha_1 - \alpha}{\alpha_1-\alpha_2}
\end{array}\right.$. 
Then, $\bP\left\{V_{n+1} \leq Q(\tilde\alpha;\hat \mF ) \right\}= \alpha$.
\end{theorem}

In this section, we presented LCP with general and potentially data-dependent $V(.) = V(.;\mZ)$, and showed that CP is its special case with $H_{ij}=1$. The discussion of this general construction is for theoretical completeness, as the general recipe described in Theorem \ref{thm:general1} or Corollary \ref{cor:generalC} is too computationally expensive: for every $Y_{n+1}=y$, we need to retrain our prediction model to get $V(.;\mZ)$.  This problem  exists in CP with data-dependent scores, and sample splitting is often used to reduce the computation cost \citep{papadopoulos2002inductive,lei2015conformal}.

For the remainder of this paper, we shift our focus to sample-splitting LCP, where we divide the observed data into a training set and calibration set. The score function $V(.)$ is estimated with the training set and considered fixed afterwards, and the PI is constructed using the fixed score function and the calibration set.

\section{Sample-splitting LCP}
\label{sec:LCPsplit}
\subsection{Sample-splitting LCP and marginal coverage guarantee}
This section considers sample-splitting LCP and develops an efficient algorithm.  In sample-splitting LCP, we divide the observed data into the training set $\mathcal{D}_0$ of size $n_0$, and calibration set $\mathcal{D}$ of size n. We first construct the score function $V(.)$ based on $\mathcal{D}_0$. For example, we may let $V(Z) = |Y - \hat \mu(X)|$ where $\hat \mu(.)$ is a prediction function for $Y$ learned using  $\mathcal{D}_0$. Since $V(Z)$ does not depend on the calibration set and the test sample, we refer to it as a data-independent score. We let $\{Z_1, \ldots, Z_n\}$  denote samples of the calibration set, and $Z_{n+1}$ the test sample.  In this setting, because $V(.)$ is fixed, the empirical distributions $\hat\mF_i$ for $i=1,\ldots, n+1$ depend on the value $y$ of a test sample $(X_{n+1},y)$ only via $v=V(X_{n+1},y)$. Thus $\tilde\alpha(y)$ as defined in Corollary \ref{cor:generalC} also depends on $y$ only via $v$. With a small abuse of notation, we will henceforth write $\tilde\alpha(v)$ in place of $\tilde{\alpha}(y)$, where $v=V(X_{n+1},y)$. To make explicit the dependence of the empirical distribution $\hat\mF_i$ on $v$, we introduce 
\begin{equation}\label{eq:hatFv}
\hat\mF_i(v) \coloneqq \hat\mF_i \mbox{ when }V_{n+1} = v.
\end{equation}
We express Theorem \ref{thm:general1} and  Corollary \ref{cor:generalC} with sample-splitting using Lemma \ref{lem:split1} below, where we can easily check that the PI for $V_{n+1}$ is an interval.
\begin{lemma}
\label{lem:split1} 
Let $V(.)$ be a fixed score function. At $V_{n+1} = v$, define $\tilde\alpha(v)$ to be the smallest value of $\tilde\alpha\in\Gamma$ such that, 
\begin{equation}
\label{eq:goal1}
\sum^{n+1}_{i=1}\frac{1}{n+1}\mathbbm{1}_{V_i \leq Q\left(\tilde\alpha; \hat\mF_i(v)\right)}  \geq \alpha.
\end{equation}
Set $C_V(X_{n+1}) = \{v: v \leq Q(\tilde\alpha(v);\hat \mF)\}$, $C(X_{n+1}) = \{y: V(X_{n+1},y)\in C_V(X_{n+1})\}$. Then $C_V(X_{n+1}) $ is an interval, and 
\[
\bP\left\{V_{n+1} \in C_V(X_{n+1}) \right\}\geq \alpha,\qquad \bP\left\{Y_{n+1} \in C(X_{n+1}) \right\}\geq \alpha.
\]
\end{lemma} 
Lemma \ref{lem:split1} is intuitively simple.  However,  even though the score function $V(.)$ is pre-specified, it is still unrealistic to compute for every possible value of $v_{n+1}=V(X_{n+1},y)$ its own value of $\tilde\alpha(v_{n+1})$.    In Section \ref{subsec:alg}, we provide an efficient implementation of LCP to tackle this problem.
\subsection{An efficient implementation of  LCP}
\label{subsec:alg}
We provide an $\mathcal{O}(n\log n)$ implementation of LCP, given pre-calculated localizer function values for each pair of calibration samples and the associated unnormalized cumulative probabilities.  

Without loss of generality, we assume  that the calibration samples are ordered by their score values and $V_1\leq V_2\leq \ldots \leq V_n$. Let $\overline{V}_i$ be the augmented observation with $\overline{V}_i = V_i$ for $i = 1,\ldots, n$, $\overline{V}_{n+1} = \infty$ and $\overline{V}_0 = -\infty$. For all $i=1,\ldots, n+1$, we define
\begin{itemize}
\item  $\ell(i)= \max\{i' \in \{1,\ldots,n\}: V_{i'}<\overline{V}_i\}$ as the largest index of values $V_{i'}$ that are smaller than $\overline{V}_i$.  In the case where all $V_i$ values are distinct, we have $\ell(i)=i-1$. We  set the maximum of an empty set as 0, so in particular, $\ell(1) = 0$ always.
\item $\theta_{i} \coloneqq \sum_{j=1}^{\ell(i)}p^{H}_{i, j}$  as the cumulative probability at the value $\overline{V}_{l(i)}$ in the distribution $\hat\mF_i(\infty)$. 
\item $\tilde \theta_{i} \coloneqq \sum_{j=1}^{\ell(i)}p_{n+1, j}^{H}$ as the cumulative probability at $\overline{V}_{l(i)}$ in the distribution $\hat\mF$. 
\item $\theta_i=\tilde\theta_{i} = 0$ if $\ell(i) = 0$. In particular, $ \theta_{1} =\tilde\theta_1 = 0$ always. 
\end{itemize}
Lemma \ref{lem:alg1} below is the foundation of our implementation.  The first part of Lemma \ref{lem:alg1} describes a formulation to construct the closure of the PI $C_V(X_{n+1})$ from Lemma \ref{lem:split1} that does not explicitly require calculation of $\tilde{\alpha}(v_{n+1})$ for different values of $v_{n+1}=V(X_{n+1},y)$. This formulation depends on a quantity $S(k)$ defined in (\ref{eq:goal1neg}). The second part of Lemma \ref{lem:alg1} gives another equivalent characterization of $S(k)$ that enables its computation for all $k=1,\ldots, n+1$ in $\mathcal{O}(n\log n)$ time.

\begin{lemma}[Practical implementation of LCP]
\label{lem:alg1}
\;
\begin{enumerate}
\item   Let $k^*$ be the largest index $k \in \{1,\ldots,n+1\}$ such that
\begin{equation}
\label{eq:goal1neg}
S(k)\coloneqq \sum_{i=1}^{n}\frac{1}{n+1}\mathbbm{1}_{V_i \leq Q(\tilde\theta_{k}; \hat{\mathcal{F}}_i(\overline{V}_{\ell(k)}))}  <\alpha.
\end{equation}
Then,  $\bar{C}_V(X_{n+1})=\{v: v\leq \overline{V}_{k^*}\}$ is the closure of $C_V(X_{n+1})$ from Lemma \ref{lem:split1}.
\item  We may partition the $n$ calibration samples into three sets: $A_1 \coloneqq \{i: p_{i,n+1}^H + \theta_i < \tilde \theta_{{i}}\}$, $A_2 \coloneqq \{i: \theta_i \geq \tilde \theta_{{i}}\}$, and $A_3 \coloneqq \{i: p_{i,n+1}^H + \theta_i \geq \tilde \theta_{{i}},\,\theta_i < \tilde \theta_{{i}}\}$. For $k = 1,\ldots, n+1$, we have
\begin{align}
\label{eq:goal1equiv}
S(k)=\sum_{i\in A_1}\frac{1}{n+1}\mathbbm{1}_{\theta_i+p_{i,n+1}^H < \tilde\theta_{k}}+\sum_{i\in A_2}\frac{1}{n+1}\mathbbm{1}_{\theta_i < \tilde\theta_{k}}+\sum_{i\in A_3}\frac{1}{n+1}\mathbbm{1}_{l(i)< {\ell(k)}}.
\end{align}
\end{enumerate}
\end{lemma}
Here, we provide some intuition for why (\ref{eq:goal1neg}) and (\ref{eq:goal1equiv}) are equivalent.  Observe that  $\tilde\theta_k$ and $\overline{V}_{\ell(k)}$ are both non-decreasing in $k$. Then the quantile $Q(\tilde\theta_k, \hat \mF_i(\overline{V}_{\ell(k)}))$ is also non-decreasing in $k$, where we recall the definition (\ref{eq:hatFv}) for $\hat \mF_i(v)$. As a result, defining the event $J_{ik} = \{V_i\leq Q(\tilde\theta_k; \hat\mF_i(\bar{V}_{\ell(k)}))\}$ in the indicator of (\ref{eq:goal1neg}), once $J_{ik}$ holds for some $k$, it holds also for all larger $k$. Thus, for each $i=1,\ldots,n$, we need only determine the smallest $k$ for which $J_{ik}$ first holds. There are two cases:
\begin{itemize}
\item If $J_{ik}$ first holds at a value $k$ with $V_i > \bar{V}_{\ell(k)}$,  by definition of $\hat\mF_i(\overline{V}_{\ell(k)})$, we need 
\[
\tilde\theta_k >\sum_{j\leq n: V_{j}< V_i}p^H_{i,j}+p^H_{i,n+1} = \theta_i+p^H_{i,n+1}.
\]
\item  If $J_{ik}$ first holds at a value $k$ with $V_i \leq \bar{V}_{\ell(k)}$, then we need instead $\tilde\theta_k >\sum_{j\leq n: V_{j}< V_i}p^H_{i,j} = \theta_i$. To guarantee that $V_i \leq \bar{V}_{\ell(k)}$, we also require $\ell(k) > \ell(i)$.
\end{itemize}
Let $k_i$ be the smallest index $k$ for which $J_{ik}$ first holds. We can show that 
\begin{itemize}
\item  $A_1$ contains all $i$ such that $V_i > \bar{V}_{\ell(k_i)}$.
\item  $A_2$ contains all $i$ such that $V_i \leq \bar{V}_{\ell(k_i)}$ and $\{\tilde\theta_{k_i} > \theta_i,\,\ell(k_i) > \ell(i)\} = \{\tilde\theta_{k_i} > \theta_i\}$.
\item  $A_3$ contains all $i$ such that $V_i \leq \bar{V}_{\ell(k_i)}$ and $\{\tilde\theta_{k_i} > \theta_i,\,\ell(k_i) > \ell(i)\} = \{\ell(k_i) > \ell(i)\}$.
\end{itemize}
This will establish the equivalence between (\ref{eq:goal1neg}) and (\ref{eq:goal1equiv}).

The desirable aspect of dividing calibration samples into $A_1,A_2,A_3$ is that we can now order the calibration samples in each set based on the values of $\theta_i+p^H_{i,n+1}$, $\theta_i$, and $l(i)$ for $A_1,A_2,A_3$ respectively, and then compute all values $S(k)$ from (\ref{eq:goal1equiv}) using a single scan through the values $k=1,\ldots,n+1$.  Algorithm \ref{alg:pathLCP} implements this idea: 
\begin{itemize}
\item Line 3 calculates $\tilde\theta_i$, $\theta_i$, and $\theta_i+p_{i,n+1}^H$ for each $i=1,\ldots,n+1$; Line 4 creates $A_1,A_2,A_3$ according to Lemma \ref{lem:alg1}.
\item Line 5 orders $i \in A_1$ by $\theta_i+p_{i,n+1}^H$, $i \in A_2$ by $\theta_i$, and $i \in A_3$ by $l(i)$. As we increase $k$, samples $i$ in each set $A_1,A_2,A_3$ will satisfy $V_i \leq Q(\tilde\theta_k; \hat\mF_i(\overline{V}_{\ell(k)}))$ sequentially.
\item Lines 7-8, 9-10 and  11-12 perform these sequential checks within each set $A_1,A_2,A_3$.
\item Finally, line 14 produces the largest $k^*$ such that (\ref{eq:goal1neg}) holds for any given target level $\alpha$.
\end{itemize}
\begin{algorithm}
    \caption{ Algorithm for LCP}
    \label{alg:pathLCP} 
    \vskip -.1in
    \hspace*{\algorithmicindent} \textbf{Input:} (1) Ordered conformal scores $V_1\leq \ldots \leq V_n$, (2) associated unnormalized cumulative probability matrix $Q_{ik} = \sum_{j=1}^k H_{ij}$ for $i, k = 1,\ldots, n$, (3) $H_{n+1, i}$ and $H_{i, n+1}$ for $i = 1,\ldots, n$, and (4) the targeted level $\alpha$.
    
    \hspace*{\algorithmicindent} \textbf{Output:} A constructed PI $C_V$ for $V_{n+1}$.
    
    $\theta_i+p_{i,n+1}^H \leftarrow\frac{Q_{i,\ell(i)}+H_{i,n+1}}{Q_{i,n}+H_{i,n+1}}$,  $\theta_i \leftarrow\frac{Q_{i,\ell(i)}}{Q_{i,n}+H_{i,n+1}}$, $\tilde\theta_i\leftarrow\frac{Q_{n+1,l(i)}}{\sum_{j=1}^{n+1} H_{n+1,j}}$ for $i=1,\ldots, n+1$.

$A_1 \leftarrow \{i: \theta_i+p_{i,n+1}^H < \tilde\theta_{i}\}$, $A_2 \leftarrow \{i: \theta_i \geq \tilde\theta_{i}\}$,  $A_3 \leftarrow \{i: \theta_i+p_{i,n+1}^H \geq \tilde\theta_{i}, \theta_i< \tilde\theta_{i}\}$.
    
     Set $\check{\theta}^{A_1},\check\theta^{A_2},\check\theta^{A_3}$ as the ordered \emph{values} of $\{\theta_i+p_{i,n+1}^H: i\in A_1\}$, $\{\theta_i:i \in A_2\}$, and $\{\ell(i): i \in A_3\}$ respectively. Set $c_m = 0$, $L_ m= |A_m|$, for $m=1, \;2,\; 3$.
    
    \For{k = 1, 2,\ldots, n, n+1}{
    
       \While{$c_1 < L_1$ and $\check{\theta}^{A_1}_{c_1+1}< \tilde\theta_k$}{
            $c_1\leftarrow c_1 + 1$;
       }

        \While{$c_2 < L_2$ and $\check{\theta}^{A_2}_{c_2+1}< \tilde\theta_k$}{
              $c_2\leftarrow c_2 + 1$;
        }

        \While{$c_3 < L_3$ and $\check{\theta}^{A_3}_{c_3+1}< {\ell(k)}$}{
           $c_3\leftarrow c_3 + 1$;
        }
         
    Set $S(k)=\frac{c_1+c_2+c_3}{n+1}$.
    }
   Set $k^* = \arg\max\{k:S(k)<\alpha\}$, and return $C_V=\{v:v\leq \overline{V}_{k^*}\}.$
\end{algorithm}

\subsection{Choice of H}
\label{subsec:method_h}
The choice of $H$ will  influence the localization. Given $d(x_1, x_2)$ as a measure of dissimilarity between two samples $x_1$, $x_2$, there are numerous ways of defining the functional form for the localizer.  In our experiments,  we consider the localizer  $H_h(x_1, x_2) =  \exp(-\frac{d(x_1, x_2)}{h})$.

A smaller $h$ results in more localization. We want to choose $h$ to have relatively narrow PIs for most samples. More specifically, we consider the following constrained objective:
\begin{align*}
J(h) & = \mbox{Average of PI}^{finite}\mbox{ length}+\lambda\times \mbox{Average of conditional  PI}^{finite}\mbox{ length's variability},\\
s.t. \;\;& \mbox{Average percent of infinite PIs is at most }\varepsilon.
\end{align*}
The parameter $\lambda$ reflects our aversion for the variability of constructed PI length at each fixed point $X_{n+1}=x$ of the sample space. We set $\lambda=1$ by default.
 
These averages are unknown and need to be estimated from the data. Recall that the score function $V(.)$ is constructed using an independent training set $\mathcal{D}_0$, whose model complexity is often tuned with cross-validation. We suggest using $\mathcal{D}_0$ and its cross-validated scores to empirically estimate the three terms in the above objective. The mathematical definitions of $J(h)$ and details of the empirical estimates are given in Appendix \ref{app:method_h}.

In low dimensions, we can have asymptotic conditional coverage as $n\rightarrow \infty$ using typical distance dissimilarities, e.g., Euclidean distance, and by choosing $h\rightarrow 0$ under suitable assumptions (see Section \ref{sec:conditionals}).  This is an ideal setting. In practice, a good user-specified dissimilarity function $d(.,.)$ will lead to improved performance in terms of constructed PI length and adaptation to the underlying heterogeneity.   Such a dissimilarity function should capture directions of feature space in which  the PI (of V) is more likely to vary.  A comprehensive and in-depth discussion of $d(.,.)$, especially in high dimensions, is beyond the scope of this paper.  In our numerical experiments,  we will define $d(.,.)$ as a weighted sum of three components: (1) $d_1(x_1, x_2) = \|\hat\rho(x_1)-\hat\rho(x_2)\|_2$, where $\hat\rho(X)$ is the estimated spread of $V(X,Y)$ conditional on $X$ \citep{lei2018distribution}; (2) $d_2(x_1, x_2) = \|P_{\parallel}(x_1 - x_2)\|_2$, where $P_{\parallel}$ is the projection onto the space spanned by the top singular vectors of the Jacobian matrix of $\hat\rho(X)$ for $X\in \mathcal{D}_0$; and (3) $d_3(x_1, x_2) = \|P_{\perp}(x_1 - x_2)\|_2$, where $P_{\perp}$ is the projection onto the space orthogonal to $P_{\parallel}$.

We include the first component since $\hat\rho(X)$ is trying to capture the heterogeneity of $V(X,Y)$. We include the second and the third components because $\hat\rho(X)$ may not fully capture this heterogeneity, so that the dissimilarity still depends on other directions of feature space. Intuitively, we can think of the projection $P_{\parallel}$ as capturing the directions of feature space in which $\rho(X)$ is more variable across the training set, and $P_{\perp}$ as capturing the remaining less important directions. We provide more details on constructing $d(.,.)$ in Appendix \ref{app:method_h}.

\section{Empirical studies: Comparison of CP and LCP}
\label{sec:empirical}
We compare LCP and CP in this section, using different numerical examples to demonstrate their differences and potential gains using LCP. We consider the usual regression problem:
\[
Y= \mu(X)+\eps, \;\eps \indep X,
\]
and four types of conformal score construction:

\noindent\textbf{Regression  score} $V^{R}(X, Y) = |Y - \hat \mu(X)|$ where $\hat \mu(X)$ is an estimate of $\mu(X)$ learned from the training set. We denote the two different procedures based on the regression score as conformalized regression (CR) and localized\&conformalized regression (LCR).

\noindent\textbf{Locally weighted regression score} $V^{R-local}(X, Y) = \frac{V^{R}(X, Y)}{\hat \rho(X)}$ ,  where $\hat\rho(X)$ is the estimated spread of $V^{R}(X,Y)$  \citep{lei2018distribution}. The locally weighted regression score also leads to two procedures: conformalized locally weighted regression (CLR) and localized\&conformalized locally weighted regression (LCLR).

\noindent\textbf{Quantile regression score} $V^{QR}(X, Y) = \max\{\hat q_{lo}(X_i) - Y, Y - \hat q_{hi}(X)\}$,  where $\hat q_{lo}(.)$ and $\hat q_{hi}(.)$ are the estimated lower and upper $\frac{\alpha}{2}$ quantiles from the training set \cite{romano2019conformalized}. The two procedures based on the quantile regression score are conformalized quantile regression (CQR), and local\&conformalized quantile regression (LCQR).

\noindent\textbf{Locally weighted quantile regression score} $V^{QR-local} =\frac{V^{QR}(X, Y)}{\hat \rho(X)}$, which combines quantile regression with the locally weighted step. The two related procedures are conformalized locally weighted quantile regression (CLQR), and local\&conformalized locally weighted quantile regression (LCLQR). 

In  Example \ref{example1}, we visually demonstrate CR and LCR to highlight the procedural differences, and compare LCP results with different values for $h$. In Example \ref{example2}, we use synthetic data and compare the performance of the eight procedures. Example \ref{example3} compares the results using four publicly available data sets from UCI.  In all empirical examples, we learn the conformal scores using a neural network with three fully connected layers and 32 hidden nodes.
\begin{example}[An illustrating example on CR and LCR] 
\label{example1}
Let $Y = \eps$, $X\sim \mathcal{N}(0, 1)$, $\eps\sim \rho(X)$, with four different cases for  $\rho(X)$: (A) $\rho(X) = \sin(X)$; (B) $\rho(X) = \cos(X)$; (C) $\rho(X) = \sqrt{|X|}$; (D) $\rho(X) = 1$. We compare CR and auto-tuned (Section \ref{subsec:method_h}) LCR, as well as results from LCR using fixed $h$ values. (The prefixed grids for $h$ can be different for different settings because they are chosen by looking at the dissimilarity measures on the training set.) The sizes for the training and calibration  sets are both 1000.  Table \ref{tab:example1a} compares CR with auto-tuned LCR, and shows the achieved coverage, percents of samples with infinite PI, and the average length of finite PI (ave.PI).  Table \ref{tab:example1b} compares LCR from using different  $h$. Figure \ref{fig:example1} provides  visual demonstrations for CR and LCR using the smallest $h$ with less than 5\% of infinite PI for LCR ($h_1$), the largest $h$ considered ($h_2$) and the auto-tuned $h$ ($h_3$). Choice of $h_1$ results in a highly localized LCR with PIs better capturing the underlying heterogeneity but potentially less stable and containing infinite PIs with higher probability, while the choice of $h_2$ results in PIs with almost no localization and almost identical to CR.
\begin{table}
\centering
\caption{Example \ref{example1}. Coverage and length comparison for CR and LCR (auto-tuned) across four simulation setups.  Column names representing LCP procedure are in bold. We also highlight the smallest ave.PI and the method that is within .05 away from it in different settings.}
\label{tab:example1a}
\begin{adjustbox}{width = .7\textwidth}
\begin{tabular}{|l|ll|ll|ll|ll|}
  \hline
 & setting A&  & setting B &  & setting C & & setting D &  \\ \hline
 & CR & \textbf{LCR} & CR &  \textbf{LCR} & CR &  \textbf{LCR} & CR &  \textbf{LCR} \\ 
  \hline
coverage & 0.95 & 0.94 & 0.95 & 0.95 & 0.95 & 0.95 & 0.94 & 0.94 \\ 
  infinite PI\% & --& 0.00 &--  & 0.00 & -- & 0.01 &--  & 0.00 \\ 
  ave.PI & 2.77 &  \textbf{2.27} & 3.14 &  \textbf{3.01} & 4.26 &  \textbf{3.15} &  \textbf{3.81} &  \textbf{3.86} \\ 
   \hline
\end{tabular}
\end{adjustbox}
\end{table}
\begin{table}
\centering
\caption{Example \ref{example1}. Comparisons of coverage, percent of infinite PI, ave.PI and ave.PI0 for different tuning parameters $h$, where ave.PI is the average length for finite PIs at the given $h$, and ave.PI0 is the average length for PIs that are finite using all $h$ considered. We highlight all ave.PI0 for $h$ no greater than the auto-tuned $\hat h$ in each setting.}
\label{tab:example1b}
\begin{adjustbox}{width=1\textwidth, height = .25\textwidth}
\begin{tabular}{|l|rrrrrrrrrrrrrrrrrrrr|}
  \hline
setting A&&&&&&&&&&&&&&&&&&&&\\  \hline
 & 0.05 & 0.07 & 0.09 & 0.13 & 0.17 & 0.22 & 0.29 & 0.39 & 0.52 & 0.69 & 0.91 & 1.21 & 1.61 & 2.14 & 2.84 & 3.78 & 5.01 & 6.66 & 8.84 & 11.74 \\ 
  \hline
coverage & 0.95 & 0.95 & 0.95 & 0.95 & 0.95 & 0.95 & 0.95 & 0.95 & 0.95 & 0.95 & 0.94 & 0.95 & 0.95 & 0.95 & 0.95 & 0.95 & 0.95 & 0.95 & 0.95 & 0.95 \\ 
  infinite PI\% & 0.23 & 0.15 & 0.07 & 0.04 & 0.02 & 0.01 & 0.01 & 0.01 & 0.00 & 0.00 & 0.00 & 0.00 & 0.00 & 0.00 & 0.00 & 0.00 & 0.00 & 0.00 & 0.00 & 0.00 \\ 
  ave.finitePI & 2.22 & 2.27 & 2.27 & 2.28 & 2.29 & 2.32 & 2.32 & 2.32 & 2.32 & 2.29 & 2.27 & 2.36 & 2.43 & 2.51 & 2.57 & 2.61 & 2.67 & 2.70 & 2.71 & 2.74 \\ 
  ave.PI0 & \textbf{2.22} &  \textbf{2.21} &  \textbf{2.20} &  \textbf{2.18} &  \textbf{2.19} &  \textbf{2.23} &  \textbf{2.24} &  \textbf{2.26} &  \textbf{2.25} &  \textbf{2.22} &  \textbf{2.22} & 2.32 & 2.40 & 2.49 & 2.55 & 2.61 & 2.66 & 2.69 & 2.70 & 2.73 \\ 
   \hline

setting B&&&&&&&&&&&&&&&&&&&&\\  \hline
 & 0.08 & 0.1 & 0.13 & 0.18 & 0.24 & 0.32 & 0.43 & 0.57 & 0.76 & 1.02 & 1.36 & 1.81 & 2.42 & 3.24 & 4.32 & 5.77 & 7.7 & 10.28 & 13.73 & 18.33 \\ 
  \hline
coverage & 0.94 & 0.95 & 0.95 & 0.94 & 0.94 & 0.95 & 0.94 & 0.94 & 0.95 & 0.95 & 0.95 & 0.95 & 0.95 & 0.95 & 0.95 & 0.95 & 0.95 & 0.95 & 0.95 & 0.95 \\ 
  infinite PI\% & 0.17 & 0.09 & 0.07 & 0.05 & 0.04 & 0.03 & 0.01 & 0.01 & 0.00 & 0.00 & 0.00 & 0.00 & 0.00 & 0.00 & 0.00 & 0.00 & 0.00 & 0.00 & 0.00 & 0.00 \\ 
  ave.finitePI & 2.81 & 2.68 & 2.67 & 2.63 & 2.65 & 2.65 & 2.66 & 2.73 & 2.82 & 2.87 & 2.98 & 3.01 & 3.03 & 3.07 & 3.08 & 3.11 & 3.12 & 3.13 & 3.13 & 3.13 \\ 
  ave.PI0&  \textbf{2.81} &  \textbf{2.81} &  \textbf{2.84} &  \textbf{2.81} &  \textbf{2.83} &  \textbf{2.83} &  \textbf{2.84} &  \textbf{2.89} &  \textbf{2.96} &  \textbf{2.98} &  \textbf{3.06} &  \textbf{3.07} & 3.07 & 3.10 & 3.11 & 3.13 & 3.13 & 3.14 & 3.13 & 3.14 \\ 
   \hline
setting C&&&&&&&&&&&&&&&&&&&&\\  \hline
 & 0.04 & 0.06 & 0.09 & 0.12 & 0.17 & 0.24 & 0.33 & 0.46 & 0.65 & 0.91 & 1.27 & 1.78 & 2.49 & 3.48 & 4.87 & 6.83 & 9.56 & 13.38 & 18.74 & 26.23 \\ 
  \hline
coverage & 0.94 & 0.94 & 0.94 & 0.94 & 0.95 & 0.95 & 0.95 & 0.95 & 0.95 & 0.95 & 0.95 & 0.95 & 0.95 & 0.95 & 0.95 & 0.95 & 0.95 & 0.95 & 0.95 & 0.95 \\ 
  infinite PI\% & 0.29 & 0.22 & 0.14 & 0.09 & 0.07 & 0.05 & 0.03 & 0.02 & 0.01 & 0.01 & 0.00 & 0.00 & 0.00 & 0.00 & 0.00 & 0.00 & 0.00 & 0.00 & 0.00 & 0.00 \\ 
  ave.finitePI & 1.88 & 2.12 & 2.42 & 2.68 & 2.81 & 2.87 & 3.02 & 3.08 & 3.15 & 3.24 & 3.34 & 3.48 & 3.62 & 3.80 & 3.94 & 4.01 & 4.10 & 4.16 & 4.19 & 4.21 \\ 
  ave.PI0 &  \textbf{1.88} &  \textbf{1.88} &  \textbf{1.87} &  \textbf{1.91} &  \textbf{1.91} &  \textbf{1.93} &  \textbf{1.99} &  \textbf{2.06} &  \textbf{2.20} & 2.44 & 2.65 & 2.97 & 3.25 & 3.57 & 3.80 & 3.90 & 4.03 & 4.12 & 4.16 & 4.19 \\ 
   \hline

setting D&&&&&&&&&&&&&&&&&&&&\\  \hline
 & 0.08 & 0.11 & 0.15 & 0.2 & 0.28 & 0.38 & 0.53 & 0.73 & 1 & 1.38 & 1.9 & 2.61 & 3.6 & 4.95 & 6.82 & 9.4 & 12.94 & 17.82 & 24.54 & 33.8 \\ 
  \hline
coverage & 0.94 & 0.95 & 0.94 & 0.94 & 0.95 & 0.94 & 0.94 & 0.94 & 0.94 & 0.94 & 0.94 & 0.94 & 0.94 & 0.94 & 0.94 & 0.94 & 0.94 & 0.94 & 0.94 & 0.94 \\ 
  infinite PI\% & 0.14 & 0.10 & 0.05 & 0.02 & 0.01 & 0.01 & 0.01 & 0.00 & 0.00 & 0.00 & 0.00 & 0.00 & 0.00 & 0.00 & 0.00 & 0.00 & 0.00 & 0.00 & 0.00 & 0.00 \\ 
  ave.finitePI & 3.79 & 3.85 & 3.92 & 3.92 & 3.92 & 3.90 & 3.89 & 3.92 & 3.89 & 3.89 & 3.89 & 3.87 & 3.86 & 3.86 & 3.88 & 3.87 & 3.87 & 3.86 & 3.85 & 3.85 \\ 
  ave.PI0 &  \textbf{3.79} &  \textbf{3.80} &  \textbf{3.80} &  \textbf{3.78} &  \textbf{3.82} &  \textbf{3.82} &  \textbf{3.82} &  \textbf{3.86} &  \textbf{3.84} &  \textbf{3.85} &  \textbf{3.86} &  \textbf{3.85} &  \textbf{3.84} &  \textbf{3.86} &  \textbf{3.87} &  \textbf{3.87} &  \textbf{3.86} &  \textbf{3.86} & 3.85 & 3.84 \\ 
   \hline
\end{tabular}
\end{adjustbox}
\end{table}
\begin{figure}
\caption{\em \em Example \ref{example1}. Confidence bands constructed using CR, LCR with different tuning parameter values for $h$  at  the targeted level $\alpha = .95$. In each of the sub-plot, we show the test data points with black dots, and the true confidence bands across different $X_{n+1}$ as the blue curves, and the estimated PIs using different methods as red curves. We sometimes encounter infinite PIs using LCP. We represent those infinite PIs by widths larger than the gray horizontal lines.}
\label{fig:example1}
\begin{center}
\includegraphics[width = .9\textwidth, height = .7\textwidth]{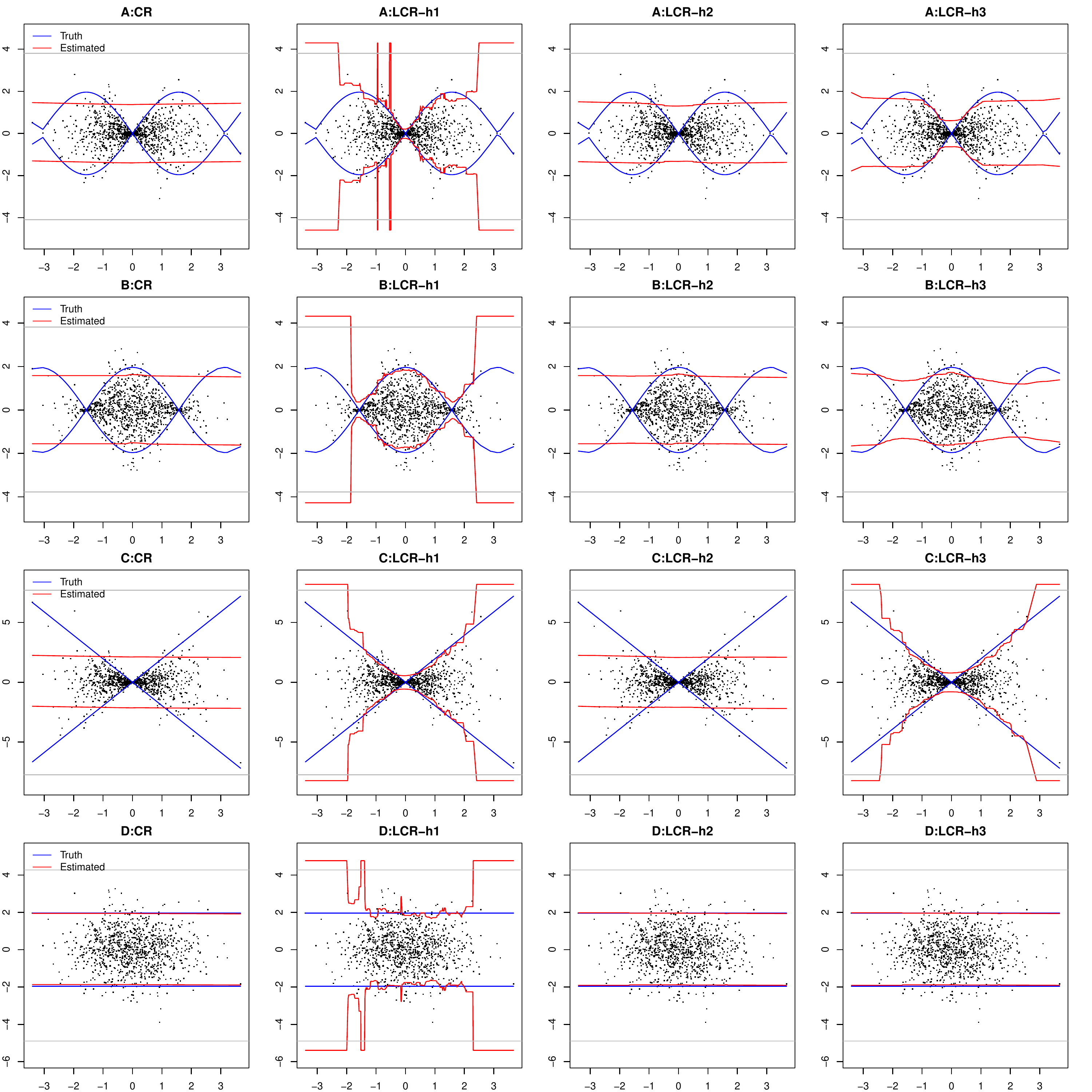} 
\end{center}
\vskip -.2in
\end{figure}
\end{example}
We do not observe an increased average PI length on samples that are well represented by the calibration set as we decrease $h$ in a wide range. A smaller $h$ makes the procedure more alert by producing PIs with infinite length for underrepresented new observations. Is this a bad thing? We believe that the answer to this question is subjective and depends on the specific task at hand.
\begin{example}[Comparisons of different procedures, synthetic data]
\label{example2}
We consider eight procedures by applying CP and LCP (auto-tuned) to four different conformal scores regarding their coverage and PI lengths at a targeted level $\alpha = 0.95$.  We consider the same simulation setup as in Example \ref{example1} except with $X\sim \rm{U}[-2,2]$.  In this example, the test observations are reasonably well represented by the calibration samples (with high probability). We do not observe samples with infinite PIs using LCP and auto-tuned $h$. Table \ref{tab:example2coverage} and \ref{tab:example2length} show the results of average coverage and average length of PI  in the four simulation settings. 
\begin{table}
\centering
\caption{Example \ref{example2}. Empirical Coverage for different procedures across different simulation settings, with a targeted level at $\alpha = 0.95$. Column names representing LCP procedure are in bold.}
\label{tab:example2coverage}
\begin{adjustbox}{height = .07\textwidth}
\begin{tabular}{lrrrrrrrr}
  \hline
 & CR & \textbf{LCR} & CLR & \textbf{LCLR} & CQR & \textbf{LCQR} & CLQR & \textbf{LCLQR} \\ 
  \hline
setting A & 0.952 & {0.955} & 0.953 & 0.954 & 0.953 & {0.955} & 0.953 & 0.954 \\ 
  setting B & 0.951 & 0.953 & 0.950 & 0.954 & 0.953 & 0.953 & 0.953 & 0.954 \\ 
  setting C & 0.948 & 0.950 & 0.949 & 0.950 & 0.950 & 0.951 & 0.951 & {0.952} \\ 
  setting D & 0.948 & 0.948 & 0.948 & {0.949} & 0.948 & {0.948} & 0.947 & 0.948 \\ 
   \hline
\end{tabular}
\end{adjustbox}
\end{table}
\begin{table}
\centering
\caption{Example \ref{example2}. Average lengths of PIs for different procedures across four different simulation settings. Column names representing LCP procedure are in bold. The smallest average PI length and those within 0.05 from it are also highlighted in each setting.}
\label{tab:example2length}
\begin{adjustbox}{width = .7\textwidth}
\begin{tabular}{lrrrrrrrr}
  \hline
 & CR &\textbf{ LCR }& CLR &\textbf{ LCLR} & CQR &\textbf{LCQR} & CLQR &\textbf{LCLQR} \\ 
  \hline
setting A & 3.27 & \textbf{2.84} & 3.05 & \textbf{2.81} & 2.87 & 2.87 & 2.88 & 2.87 \\ 
  setting B & 2.86 &\textbf{2.19} & 2.54 & \textbf{2.20} & 2.26 & 2.27 & 2.27 & 2.27 \\ 
  setting C & 4.95 & \textbf{3.91} & 4.42 & \textbf{3.90 }& \textbf{3.94} & \textbf{3.95} & 4.03 & 4.02 \\ 
  setting D & \textbf{3.88} &  \textbf{3.90} &  \textbf{3.89} &  \textbf{3.90} &  \textbf{3.92} &  \textbf{3.93} &  \textbf{3.92} &  \textbf{3.93} \\ 
   \hline
\end{tabular}
\end{adjustbox}
\end{table}
\end{example}
\begin{example}[Performance comparison on four UCI datasets]
\label{example3}
We investigate the performances of eight procedures (with auto-tuned LCP) on four UCI datasets \citep{rana2013physicochemical}:  CASP \citep{yeh1998modeling}, Concrete \citep{yeh1998modeling}, Facebook variant 1 (facebook1) and variant 2 (facebook2)\citep{singh2015comment, singh2015facebook}. The sizes of samples and features are (45730,9), (1030,8), (40949, 53), (81312, 53) for the four datasets respectively.

We subsample 5000  training/calibration samples without replacement from CASP, Facebook variant 1, and Facebook variant 2,  and  400 training/calibration samples from the Concrete dataset. We construct PIs using the remaining samples for each data set and repeat it 20 times. Tables \ref{tab:example3coverage} - \ref{tab:example3length} show the results of average coverage and average length of finite PI for the four data sets. The percent of infinite PIs ranges from 0\% to 3\% for different LCP constructions. The samples with infinite PI using LCP methods on the Facebook datasets tend to have wider PIs, hence,  we also show the average PI length using only samples with finite PI  from all procedures for a fair comparison.  
\begin{table}
\centering
\caption{Example \ref{example3}. Empirical Coverage for different procedures, with a targeted level at $\alpha = 0.95$. Column names representing LCP procedure are highlighted.}
\label{tab:example3coverage}
\begin{adjustbox}{height = .07\textwidth}
\begin{tabular}{lrrrrrrrr}
  \hline
 & CR & \textbf{LCR} & CLR & \textbf{LCLR} & CQR & \textbf{LCQR} & CLQR &\textbf{LCLQR} \\ 
  \hline
CASP & 0.949 & 0.950 & 0.950 & 0.950 & 0.950 & 0.950 & 0.950 & 0.950 \\ 
  Concrete & 0.947 & 0.949 & 0.943 & 0.947 & 0.951 & 0.953 & 0.952 & 0.954 \\ 
  facebook1 & 0.949 & 0.950 & 0.949 & 0.949 & 0.950 & 0.951 & 0.950 & 0.951 \\ 
  facebook2 & 0.951 & 0.951 & 0.951 & 0.951 & 0.953 & 0.952 & 0.953 & 0.952 \\ 
   \hline
\end{tabular}
\end{adjustbox}
\end{table}
\begin{table}
\centering
\caption{Example \ref{example3}.  The upper half shows the average length using samples with finite PIs for the given procedure. The lower half shows the average PI length on the common set of samples with finite PI for all procedures. Column names representing LCP procedure are in bold. The smallest average PI length and those within 0.05 from it are also highlighted in each setting.}
\label{tab:example3length}
\begin{adjustbox}{width = .8\textwidth}
\begin{tabular}{lrrrrrrrr}
\hline
Procedure-specific samples\\
  \hline
 & CR & \textbf{LCR} & CLR & \textbf{LCLR} & CQR & \textbf{LCQR} & CQLR & \textbf{LCQLR} \\ 
  \hline
CASP & 3.03 & 2.82 & 2.88 & 2.81 & \textbf{2.65} & \textbf{2.64} & \textbf{2.65} & \textbf{2.64} \\ 
  Concrete & 1.52 & \textbf{1.48} & \textbf{1.45} & \textbf{1.46} &1.58 & {1.59} & 1.58 & 1.57 \\ 
  facebook1 & 1.12 & \textbf{0.67} & $>100$ & \textbf{0.69} & 0.92 & 0.92 & 0.92 & 0.92 \\ 
  facebook2 & 1.05 & \textbf{0.65} & $>100$ & \textbf{0.65} & 0.92 & 0.92 & {0.93} & 0.92 \\ 
   \hline
Common samples\\
  \hline
 & CR &  \textbf{LCR} & CLR &  \textbf{LCLR} & CQR &  \textbf{LCQR} & CQLR &  \textbf{LCQLR} \\ 
  \hline
CASP & 3.03 & 2.82 & 2.88 & 2.81 & \textbf{2.65} & \textbf{2.64} & \textbf{2.65} &\textbf{2.64} \\ 
  Concrete & 1.52 & \textbf{1.47} & \textbf{1.43}& \textbf{1.45} & 1.56 & { 1.56} & 1.56 & 1.55 \\ 
  facebook1 & 1.12 & \textbf{0.66} & 0.95 &\textbf{0.66} & 0.74 & 0.74 & 0.74 & 0.74 \\ 
  facebook2 & 1.05 & \textbf{0.62} & 0.87 & \textbf{0.63} & 0.74 & 0.74 & 0.74 & 0.74 \\ 
   \hline
\end{tabular}
\end{adjustbox}
\end{table}
The CQR procedure has been shown as a top performer in \cite{romano2019conformalized} and \cite{sesia2020comparison}.  Our numerical experiments also confirm that it has an overall better performance than the CLR. Not only is LCP framework conceptually novel; it uses the estimated spread $\hat \rho(X)$ in a more robust way. When combined with $V^{R}$ and $V^{R-local}$, the average PI lengths are smaller for three out of the four real data sets compared with the CLR procedure. In particular, LCR and LCLR are even noticeably better than CQR in the two Facebook examples. 
\end{example}

\section{Local behavior of LCP}
\label{sec:conditionals}
In this section, we consider asymptotic and approximate conditional coverage properties for LCP, as well as for a simplified version of LCP that uses the choice $\tilde\alpha=\alpha$. We have shown in Proposition \ref{exm:counter2} that choosing $\tilde\alpha=\alpha$ does not yield a distribution-free coverage guarantee. Our results here indicate that this choice may lead to asymptotic or approximate conditional coverage, under certain assumptions.

For simplicity, in this section we restrict attention to a localizer $H(x_1, x_2) =\exp(-\frac{d(x_1, x_2)}{h_n})$ where  $h_n$ is an $n$-dependent bandwidth parameter, and $d(x_1, x_2)\geq 0$ is a measure of dissimilarity satisfying $d(x_1, x_1) = 0$.\\

{\bf Asymptotic conditional coverage.}
Non-trivial finite sample and distribution-free conditional coverage is impossible for continuous distributions \citep{lei2014distribution,vovk2012conditional}. Thus, it is common to consider asymptotic conditional coverage under proper assumptions on $\mathcal{P}_{XY}$.  Different conformal score constructions with such asymptotic conditional coverage are studied in the literature. For instance, in \cite{izbicki2019flexible}, the authors consider using the estimated conditional density as the conformal score, and in \cite{romano2019conformalized}, the authors use the conformal score based on estimated quantile functions. Here, we consider the asymptotic behavior of LCP.

\begin{assumption}
\label{ass:regularity1}
$X$ has continuous distribution on $[0,1]^p$, and $V(Z)$ has continuous distribution conditional on $X=x$. Furthermore, there exist constants $L>0$ and $\beta\geq 0$ such that the density of $X$ satisfies $p_X(x)\geq 1\slash L$ for all $x \in [0,1]^p$, and
\begin{itemize}
\item[(i)]  The conditional distribution of $V$ given $X$ satisfies for all $x,x' \in [0,1]^p$
\[
\max_{v\in \real}\left|P_{V|x}(v)-P_{V|x'}(v)\right|\leq L d(x, x'),
\]
where $P_{V|x}(v)$ is the probability that $V(Z)\leq v$ conditional on $X= x$.
\item[(ii)]  $\bP(X\in \{x:d( x_0, x) \leq \varepsilon\}) \geq \varepsilon^{\beta}/L$ for all $\varepsilon \leq h_n$ and all $x_0\in [0,1]^p$.

\item[(iii)] $h_n$ is chosen such that $h_n\rightarrow 0$ and $(n h_n^\beta\slash \ln n)\rightarrow  \infty$ as $n\rightarrow \infty$.
\end{itemize}
\end{assumption}
Under this assumption, statement (\ref{eq:thmAsy1}) of the following theorem guarantees that LCP with $\tilde\alpha(v)$ chosen as in Lemma \ref{lem:split1} achieves asymptotic conditional coverage at the target level $\alpha$. Furthermore, statements (\ref{eq:thmAsy2}) and (\ref{eq:thmAsy3}) show that $\tilde\alpha(v)$ converges to $\alpha$ in probability asymptotically, and asymptotic conditional coverage holds also if LCP is applied with the simpler choice $\tilde\alpha=\alpha$. 
\begin{theorem}
\label{thm:AsymptoticConditional}
Define $\tilde\alpha(v)$ and $C_V(X_{n+1})$  as in Lemma \ref{lem:split1}. Under Assumption \ref{ass:regularity1}, for any $x_0 \in [0,1]^p$, we have
\begin{align}
&\lim_{n\rightarrow\infty}\bP(V_{n+1}\in C_V(X_{n+1})|X_{n+1}=x_0) =\alpha,\label{eq:thmAsy1}\\
&\lim_{n\rightarrow\infty}\bP(V_{n+1}\leq Q(\alpha; \hat\mF)|X_{n+1}=x_0) = \alpha,\label{eq:thmAsy2}\\
&\lim_{n\rightarrow\infty}\bP(\max_{v}|\tilde\alpha(v)-\alpha|<\varepsilon|X_{n+1}=x_0) =1, \mbox{ for all }\varepsilon > 0.\label{eq:thmAsy3}
\end{align}
\end{theorem}

In Assumption \ref{ass:regularity1}, the measure $d(x, x')$ can be defined to capture the directions where the conditional distribution of $V$ given $X$ is more likely to change as we vary $X$. Assumption \ref{ass:regularity1} (i) allows more variability in some directions and less in others based on how the data is generated, and scales better with the dimension compared to a symmetric distance such as the Euclidean distance. Assumption \ref{ass:regularity1} (ii)  assumes that $d(x_0,x)$ has enough concentration around 0, and it holds for a typical dissimilarity measure in low dimensions. In high dimensions, this assumption holds if $d(\cdot,\cdot)$ emphasizes a few directions instead of treating all directions equally. For example, if $d(x_0,x) = |x_j - x_{0,j}|$ depends only on feature $j$, then $\bP(\{X: d(x_0,X)\leq \varepsilon\})\geq \frac{1}{L}\varepsilon$ for some large constant $L$. Assumption \ref{ass:regularity1} (iii) requires $h_n$ to decay to 0 at a sufficiently slow rate. This is so that, combined with Assumption \ref{ass:regularity1} (ii), we may ensure that $\sum_{i=1}^{n} H(X_j, X_i) \to \infty$ for all $j=1,\ldots,n+1$, with high probability. In particular, a setting such as described in Proposition \ref{exm:counter1} cannot occur.\\
 
{\bf Approximate conditional coverage.}
In  \cite{vovk2012conditional} and \cite{lei2014distribution}, the authors  partition the feature space into  $K$ finite subsets and apply conformal inference to each of the subsets: This guarantees $\bP\{Y_{n+1}\in \hC(X_{n+1})|X_{n+1}\in \mathcal{X}_k\} \geq \alpha$ for all $k = 1,2,\ldots, K$ and some fixed partition $\cup^K_{k=1} \mathcal{X}_k = \real^p$.   In \cite{barber2019limits}, the authors consider a potentially stronger version where different regions $\mathcal{X}_k$ may overlap.  \cite{barber2019conformal} introduce a different notion of approximate conditional coverage, where instead of finding $C(x_0)$ that achieves conditional coverage of $Y_{n+1}$ given $X_{n+1}=x_0$, the authors consider $C(x_0)$ that covers $\tilde Y$ whose feature value $\tilde X$ is distributed according to some locally weighted distribution around  $x_0$. (See Eqs.\ (18--19) of \cite{barber2019conformal}.) When this weighted distribution becomes increasingly concentrated around $x_0$, the distribution of $\tilde Y$ intuitively approaches the conditional distribution of $Y_{n+1}$, so this serves as an approximation to conditional coverage. Here, we show that for a local weighting given by $H(x_0,x)=\exp(-\frac{d(x_0, x)}{h_n})$, an LCP procedure using this same $H(\cdot,\cdot)$ as its localizer and $\tilde\alpha =\alpha$ can achieve this guarantee for every fixed $x_0$.

\begin{theorem}
\label{thm:ApproximateConditional}
Fix any $x_0$. Define the weighted distribution $\frac{d \tilde \mP^{x_0}_X(x)}{d x} \propto \frac{d \mP_X(x)}{d x}H(x_0, x)$.
Conditional on $X_{n+1}=x_0$, let $\tilde Z=(\tilde X, \tilde Y)\sim \tilde \mP_{X}^{x_0}\times \mP_{Y|X}$. 
Define
\begin{align*}
\varepsilon(X_{n+1})&=\max_{x:d(X_{n+1},x)<\infty}\max_y|V(X_{n+1},y)-V(x,y)|,\\
\tilde C(X_{n+1}) &= \left\{y:V(X_{n+1},y)\leq Q(\alpha; \hat\mF)+\varepsilon(X_{n+1})\right\},
\end{align*}
where $d(.,.)$ is the dissimilarity measure that defines $H(.,.)$. Then
 \begin{equation}
\label{eq:approxmateLCP}
\bP\{V(\tilde Z)\leq Q(\alpha; \hat\mF)|X_{n+1} = x_0\}\geq \alpha, \qquad\; \bP\{\tilde Y\in \tilde C(X_{n+1})|X_{n+1} = x_0\}\geq \alpha.
\end{equation}
\end{theorem}
The interval $\tilde{C}(X_{n+1})$ above remains a PI at $X_{n+1}$, and does depend on $\tilde{X}$ in its construction. The term $\varepsilon(X_{n+1})$ is introduced to bound the discrepancy in the score function as we vary $x$ in a defined neighborhood around $X_{n+1}=x_0$ with $d(x_0, x)<\infty$. The value of $\varepsilon(X_{n+1})$ depends only on the score function $V(.)$ and our definition of $d(.,.)$, not the data distribution $\mP_{XY}$.  For example, we can choose $d(.,.)$ to exclude samples that are far from each other by setting $d(x_0, x) = 0$ when $\|x_0 - x\|_2\leq h$ and $d(x_0, x) = \infty$ otherwise.  In this case, when $V(Z) = |\mu(X) - Y|$ is the regression score, we have $\varepsilon(X_{n+1})\leq \max_{\|X_{n+1}-x\|_2\leq h}|\mu(X_{n+1})-\mu(x)|$ by triangle inequality.

\section{Discussion}
We propose LCP as an extension to the conventional conformal prediction framework, which uses a weighted empirical distribution around the test sample. In our numerical experiments, LCP improves over CP when there is heterogeneity in the distribution of the score function, and the localizer in LCP is defined by a dissimilarity measure $d(.,.)$ that captures the relevant directions of such heterogeneity.  Otherwise, auto-tuned LCP ends up with PIs very similar to those from CP, given the same conformal score function. Thus, ignoring the computational cost, there is little loss in replacing CP with LCP.

One  downside of LCP is its computation compared to CP. The bulk of the additional computation lies in calculating and sorting the weights for the empirical distributions. One future direction is to reduce the computational cost of LCP  for a huge calibration set. For example, we may combine LCP with proper clustering methods or estimate an approximated cumulative probability matrix using machine learning methods to reduce the computational cost. 

CP has been used in classification problems for outlier detection \citep{hechtlinger2018cautious,guan2019prediction}.   LCP may also be a useful framework for making predictions in the presence of outliers. When choosing a suitably small $h$, LCP becomes sensitive to outliers while not increasing much the length of PI for  test samples well-represented by calibration data. 

In this paper, we considered the one-dimensional regression response. CP has also been applied to other data types, including survival data and data with multi-dimensional responses \citep{candes2021conformalized, izbicki2019flexible, feldman2021calibrated}. For multi-dimensional responses, a rectangular region formed by outer products of PIs of the individual responses does not capture potential relationships between different responses.  Various authors have worked on constructions of PIs for multi-dimensional responses to address this issue \citep{paindaveine2011directional, kong2012quantile} and \cite{feldman2021calibrated} has incorporated such constructions into CP recently. Another direction for future work is to apply the idea of LCP in similar contexts.

\bibliographystyle{agsm}
\bibliography{distributionShift}

@article{izbicki2019flexible,
  title={Flexible distribution-free conditional predictive bands using density estimators},
  author={Izbicki, Rafael and Shimizu, Gilson T and Stern, Rafael B},
  journal={arXiv preprint arXiv:1910.05575},
  year={2019}
}

@article{kong2012quantile,
  title={Quantile tomography: using quantiles with multivariate data},
  author={Kong, Linglong and Mizera, Ivan},
  journal={Statistica Sinica},
  pages={1589--1610},
  year={2012},
  publisher={JSTOR}
}

@article{paindaveine2011directional,
  title={On directional multiple-output quantile regression},
  author={Paindaveine, Davy and {\v{S}}iman, Miroslav},
  journal={Journal of Multivariate Analysis},
  volume={102},
  number={2},
  pages={193--212},
  year={2011},
  publisher={Elsevier}
}

@article{gupta2021nested,
  title={Nested conformal prediction and quantile out-of-bag ensemble methods},
  author={Gupta, Chirag and Kuchibhotla, Arun K and Ramdas, Aaditya},
  journal={Pattern Recognition},
  pages={108496},
  year={2021},
  publisher={Elsevier}
}

@article{feldman2021calibrated,
  title={Calibrated Multiple-Output Quantile Regression with Representation Learning},
  author={Feldman, Shai and Bates, Stephen and Romano, Yaniv},
  journal={arXiv preprint arXiv:2110.00816},
  year={2021}
}

@article{candes2021conformalized,
  title={Conformalized survival analysis},
  author={Cand{\`e}s, Emmanuel J and Lei, Lihua and Ren, Zhimei},
  journal={arXiv preprint arXiv:2103.09763},
  year={2021}
}

@article{sesia2020comparison,
  title={A comparison of some conformal quantile regression methods},
  author={Sesia, Matteo and Cand{\`e}s, Emmanuel J},
  journal={Stat},
  volume={9},
  number={1},
  pages={e261},
  year={2020},
  publisher={Wiley Online Library}
}

@article{rana2013physicochemical,
  title={Physicochemical properties of protein tertiary structure data set},
  author={Rana, PS},
  journal={UCI Machine Learning Repository},
  year={2013}
}

@article{yeh1998modeling,
  title={Modeling of strength of high-performance concrete using artificial neural networks},
  author={Yeh, I-C},
  journal={Cement and Concrete research},
  volume={28},
  number={12},
  pages={1797--1808},
  year={1998},
  publisher={Elsevier}
}

@inproceedings{singh2015comment,
  title={Comment volume prediction using neural networks and decision trees},
  author={Singh, Kamaljot and Sandhu, Ranjeet Kaur and Kumar, Dinesh},
  booktitle={IEEE UKSim-AMSS 17th International Conference on Computer Modelling and Simulation, UKSim2015 (UKSim2015)},
  year={2015}
}

@article{singh2015facebook,
  title={Facebook comment volume prediction},
  author={Singh, Kamaljot},
  journal={International Journal of Simulation: Systems, Science and Technologies},
  volume={16},
  number={5},
  pages={16--1},
  year={2015}
}

@article{lei2018distribution,
  title={Distribution-free predictive inference for regression},
  author={Lei, Jing and G'Sell, Max and Rinaldo, Alessandro and Tibshirani, Ryan J and Wasserman, Larry},
  journal={Journal of the American Statistical Association},
  volume={113},
  number={523},
  pages={1094--1111},
  year={2018},
  publisher={Taylor \& Francis}
}

@inproceedings{romano2019conformalized,
  title={Conformalized quantile regression},
  author={Romano, Yaniv and Patterson, Evan and Candes, Emmanuel},
  booktitle={Advances in Neural Information Processing Systems},
  pages={3538--3548},
  year={2019}
}

@article{hechtlinger2018cautious,
  title={Cautious deep learning},
  author={Hechtlinger, Yotam and P{\'o}czos, Barnab{\'a}s and Wasserman, Larry},
  journal={arXiv preprint arXiv:1805.09460},
  year={2018}
}

@article{guan2019prediction,
  title={Prediction and outlier detection in classification problems},
  author={Guan, Leying and Tibshirani, Rob},
  journal={arXiv preprint arXiv:1905.04396},
  year={2019}
}

@article{lei2015conformal,
  title={A conformal prediction approach to explore functional data},
  author={Lei, Jing and Rinaldo, Alessandro and Wasserman, Larry},
  journal={Annals of Mathematics and Artificial Intelligence},
  volume={74},
  number={1-2},
  pages={29--43},
  year={2015},
  publisher={Springer}
}

@inproceedings{papadopoulos2002inductive,
  title={Inductive confidence machines for regression},
  author={Papadopoulos, Harris and Proedrou, Kostas and Vovk, Volodya and Gammerman, Alex},
  booktitle={European Conference on Machine Learning},
  pages={345--356},
  year={2002},
  organization={Springer}
}

@article{barber2019limits,
  title={The limits of distribution-free conditional predictive inference},
  author={Barber, Rina Foygel and Candes, Emmanuel J and Ramdas, Aaditya and Tibshirani, Ryan J},
  journal={arXiv preprint arXiv:1903.04684},
  year={2019}
}

@inproceedings{vovk2012conditional,
  title={Conditional validity of inductive conformal predictors},
  author={Vovk, Vladimir},
  booktitle={Asian conference on machine learning},
  pages={475--490},
  year={2012}
}

@article{lei2014distribution,
  title={Distribution-free prediction bands for non-parametric regression},
  author={Lei, Jing and Wasserman, Larry},
  journal={Journal of the Royal Statistical Society: Series B (Statistical Methodology)},
  volume={76},
  number={1},
  pages={71--96},
  year={2014},
  publisher={Wiley Online Library}
}

@article{barber2019conformal,
  title={Conformal Prediction Under Covariate Shift},
  author={Barber, Rina Foygel and Candes, Emmanuel J and Ramdas, Aaditya and Tibshirani, Ryan J},
  journal={arXiv preprint arXiv:1904.06019},
  year={2019}
}

@article{shafer2008tutorial,
  title={A tutorial on conformal prediction},
  author={Shafer, Glenn and Vovk, Vladimir},
  journal={Journal of Machine Learning Research},
  volume={9},
  number={Mar},
  pages={371--421},
  year={2008}
}

@book{vovk2005algorithmic,
  title={Algorithmic learning in a random world},
  author={Vovk, Vladimir and Gammerman, Alex and Shafer, Glenn},
  year={2005},
  publisher={Springer Science \& Business Media}
}

@article{vovk2009line,
  title={On-line predictive linear regression},
  author={Vovk, Vladimir and Nouretdinov, Ilia and Gammerman, Alex and others},
  journal={The Annals of Statistics},
  volume={37},
  number={3},
  pages={1566--1590},
  year={2009},
  publisher={Institute of Mathematical Statistics}
}

\newpage
\appendix

In this Supplement, we describe a few supplemental Lemmas used in our proofs to results in the main paper in Appendix \ref{app:supp_lemma}. We then give proofs to results in the main paper in Appendix \ref{app:proofmain} and proofs to the supplemental Lemmas in Appendix \ref{app:proofapp}. We provide details of the construction of the dissimilarity measure $d(.,.)$ used in this paper and the automatic choice of $h$ in Appendix \ref{app:method_h}.

Without loss of generality, we always assume that $V_1\leq V_2\leq \ldots\leq V_{n}$ in this supplement.
\section{A collection of supplemental Lemmas}
\label{app:supp_lemma}
Lemma \ref{lem:infequiv} describes the elementary relationship  used in the proof from previous work on weighted conformal prediction \citep{barber2019conformal}, and we state it here for the reader's convenience.   Lemma \ref{lem:monotone} states the monotone dependence of $Q(\tilde\alpha;\hat\mF_i(v))$ on $\tilde\alpha$ or $v$. Lemma \ref{lem:conditional} is a core Lemma on the marginal coverage guarantee for LCP  with strategically chosen $\tilde\alpha$. Lemma \ref{lem:alpha_bound} collects basic bounds used in the proofs of Theorem  \ref{thm:AsymptoticConditional}. 

\begin{lemma}[A.1]
\label{lem:infequiv}
For any $\alpha$ and sequence $\{V_1,\ldots, V_{n+1}\}$, we have
\[
V_{n+1} \leq Q(\alpha; \sum^n_{i=1} w_i \delta_{V_i} + w_{n+1}\delta_{V_{n+1}}) \Leftrightarrow V_{n+1} \leq Q(\alpha; \sum^n_{i=1} w_i \delta_{V_i} + w_{n+1}\delta_\infty),
\]
where $\sum^n_{i=1} w_i \delta_{V_i} + w_{n+1}\delta_{V_{n+1}}$ and $\sum^n_{i=1} w_i \delta_{V_i} + w_{n+1}\delta_\infty$ are some weighted empirical distributions with weights $w_i\geq 0$ and $\sum^{n+1}_{i=1} w_i = 1$.

\end{lemma}

\begin{lemma}[A.2]
\label{lem:monotone}
Suppose $\{V_i, i=1,\ldots, n\}$, the target level $\alpha$, and empirical weights $p_{ij}^{H}$ are given. Then,

(i) Given $V_{n+1}$, $Q\left(\tilde\alpha;\hat\mF_i(V_{n+1})\right)$ for $i=1,\ldots,n+1$ and $Q\left(\tilde\alpha;\hat\mF\right)$ are non-decreasing, right-continuous and piece-wise constant on $\tilde\alpha$, and with value changing only at  the cumulative probabilities at different $V_i$.

(ii) Given $\tilde\alpha$, $Q\left(\tilde\alpha;\hat\mF_i(v)\right)$ is non-decreasing on $v$ for $i=1,\ldots, n+1$.

(iii) If $V_{n+1} = v$ is accepted in the $C_V(X_{n+1})$ in Lemma \ref{lem:split1}, then $v'$ is accepted for any $v' \leq v$.

\end{lemma}

\begin{lemma}[A.3]
\label{lem:conditional}
Let $V_i = V(Z_i; \mZ)$ be the score for sample $i$, and $Z_i$ is i.i.d generated for $i=1,\ldots,n+1$. For any event 
\[
\mathcal{T} \coloneqq \left\{\{Z_i, i = 1,\ldots, n+1\}= \{z_i\coloneqq (x_i, y_i), i = 1,\ldots, n+1\}\right\},
\]
we have
\[
\bP\{V_{n+1} \leq Q(\tilde\alpha; \sum^{n+1}_{i=1} p^{H}_{n+1,i} \delta_{V_{i}}) |\mathcal{T}\}= \bE\left\{\frac{1}{n+1}\sum^{n+1}_{i=1}\mathbbm{1}_{v_i \leq v^*_i}|\mathcal{T}\right\},
\]
where $v_i = V(z_i; (z_1,\ldots, z_n, z_{n+1}))$,  $v^*_{i} =  Q(\tilde\alpha; \sum^{n+1}_{j=1}p^{H}_{i,j} \delta_{v_j})$ for $i = 1,2,\ldots, n+1$, and $\tilde\alpha$ can be random but is independent of the data conditional on $\mathcal{T}$. The expectation on the right side  is taken over the randomness of $\tilde\alpha$ conditional on $\mathcal{T}$.

\end{lemma}

\begin{lemma}[A.4]
\label{lem:alpha_bound}
Suppose that Assumption \ref{ass:regularity1} holds and $V(.)$ is a fixed function.   For any $x_0$, define {\small$B(x_0) = \sum_{j=1}^{n+1} H(x_0, X_i)$, $ \Delta(x_0, X)= H(x_0, X)\max_{v}|P_{V|X}( v) - P_{V|x_0}(v)|$} and {\small $\Delta(x_0)=\sum_{i=1}^n\Delta(x_0, X_i)$}. Then,
\begin{enumerate}
\item[(i)]  There exists a constant $C>0$ such that, for all $x_0\in [0,1]^p$, we have
\[
\bP(B(x_0)\leq \frac{nh_n^{\beta}}{2eL})\leq \exp(-\frac{n h_n^\beta}{8L}),\qquad \frac{\Delta(x_0)}{B(x_0)\vee (nh_n^\beta)}\leq C h_n \ln (h_n^{-1}).
\]
\item[(ii)] Set $B_i = B(X_i)$, $R_i = \frac{\sum_{j\neq i}\left(\mathbbm{1}_{V_j < V_i} - P_{V|X_j}(V_i)\right)}{B_i}$. Then, for all $V_i$ and $i=1,\ldots, n+1$, we have
\[
\bP( |R_i|\geq \sqrt{\frac{\ln n}{B_i}} |\mX, V_i)\leq \frac{2}{n^2}.
\]
\end{enumerate}
\end{lemma}

\section{Proofs Propositions, Lemmas and Theorems}
\label{app:proofmain}
In this section, we provide proofs omitted from the main paper.  We first give arguments to Proposition \ref{exm:counter1} and Proposition \ref{exm:counter2} for the counterexamples. We then present proofs to Theorem \ref{thm:general1}, Theorem \ref{thm:general2}, Lemma \ref{lem:split1}, Lemma \ref{lem:alg1} that characterize the marginal behavior of LCP and our implementation. After that, we prove Theorem \ref{thm:AsymptoticConditional} - \ref{thm:ApproximateConditional} on the asymptotic and local behaviors of LCP-type procedures.
\section*{Proofs of the counter examples}
\subsection{Proof of Proposition \ref{exm:counter1}}
\begin{proof}
When $\sum_{i=1}^{n+1} H(X_{n+1}, X_i)< \frac{1}{1-\alpha}$, by definition, we have
\[
\sum_{i=1}^{n} p_{n+1,i}^H=\frac{\sum_{i=1}^n H_{n+1, i}}{\sum_{i=1}^{n+1}H_{n+1,i}}< \frac{\frac{1}{1-\alpha}-1}{\frac{1}{1-\alpha}}=\alpha.
\]
We thus have $Q(\alpha; \hat\mF) = \infty$, and consequently, 
\begin{align*}
&\bP(Q(\alpha; \hat\mF) = \infty) = \bP(\sum_{i=1}^{n} p_{n+1,i}^H<\alpha) =  \bP(\sum_{i=1}^{n+1} H(X_{n+1}, X_i)<\frac{1}{1-\alpha}) \geq \varepsilon.\\
& \bP(Y_{n+1}\in C(X_{n+1}))\geq \bP(Q(\alpha; \hat\mF) = \infty) \geq \varepsilon.
\end{align*}
\end{proof}
\subsection{Proof of Proposition \ref{exm:counter2}}
\begin{proof}
For $X_{n+1}\in \{\pm e_j, j = 1,\ldots, p\}$,  let  $n_0$ is the number of samples with $X_i = 0$ and $n_1$ is the number of samples with $X_i = X_{n+1}$.  The achieved conditional coverage at $\tilde\alpha = \alpha$ given $\mX=X_{1:(n+1)}$ can be upper bounded as below:
\begin{align}
\label{eq:eq_counter2_1}
\bP(V_{n+1}\leq Q(\alpha; \hat \mF)|\mX)& = \bP\left(V_{n+1}\leq Q(\alpha; \frac{1}{n_1 + n_0+1}\sum_{i: X_i\in \{0, X_{n+1}\}}\delta_{V_i}+\frac{1}{n_1+n_0+1}\delta_{\infty})\bigg\rvert\mX\right)\notag\\
& \overset{(a)}{\leq} \frac{1}{n_1+1}+\frac{n_0+n_1+1}{n_1+1}[\alpha - \frac{n_0}{n_0+n_1+1}]_+,
\end{align}
Step (a) holds because:
\begin{itemize}
\item When $\alpha\leq \frac{n_0}{n_1+n_0+1}$, the $\alpha$ quantile of the weighted empirical distribution is 0, and we will have 0 coverage for $X_{n+1}\neq 0$ and (\ref{eq:eq_counter2_1}) is true.
\item  When $\alpha>\frac{n_0}{n_1+n_0+1}$, the $\alpha$ quantile of the weighted empirical distribution in  (\ref{eq:eq_counter2_1})  is the $\lceil (n_1+n_0+1)\alpha\rceil - n_0$ largest value in $\{V_{i}:X_i = X_{n+1}\}\cup V_{\infty}$, which is the $\frac{\lceil (n_1+n_0+1)\alpha\rceil - n_0}{n_1+1}$ quantile of the unweighted empirical distribution formed by $\{V_{i}:X_i = X_{n+1}\}\cup V_{\infty}$. By Lemma \ref{lem:infequiv}, $\{V_{n+1}\leq Q(t; \{V_{i}:X_i = X_{n+1}\}\cup V_{\infty})\}\Leftrightarrow \{V_{n+1}\leq Q(t; \{V_{i}:X_i = X_{n+1}\}\cup V_{n+1})\}$. Hence,  we have
{\small
\begin{align}
\bP(V_{n+1}\leq Q(\alpha; \hat \mF)|\mX)&= \bP\left(V_{n+1}\leq Q(\frac{\lceil (n_1+n_0+1)\alpha\rceil - n_0}{n_1+1}; \{V_{i}:X_i = X_{n+1}\}\cup V_{n+1})\}\bigg\rvert\mX\right)\notag\\
&\overset{(b)}{=}\frac{\lceil (n_1+n_0+1)\alpha\rceil - n_0}{n_1+1}\leq \frac{1}{n_1+1}+\frac{n_0+n_1+1}{n_1+1}(\alpha-\frac{n_0}{n_0+n_1+1}),\notag
\end{align}
}
where step $(b)$ uses the fact that $V_i\sim \operatorname{Unif}[-1,1]$ for all $i$ with $X_i = X_{n+1}$. Hence,  (\ref{eq:eq_counter2_1}) holds.
\end{itemize}
Next, we marginalize over $X_{1:n}$ but conditional on $m = n_0+n_1$ (the total number of samples with $X_i\in \{0, X_{n+1}\}$). From (\ref{eq:eq_counter2_1}):
\begin{align}
\label{eq:counter2_eq1}
\bP(V_{n+1}\leq Q(\alpha;\hat\mF)|m, X_{n+1})&\leq \bE[\frac{1}{n_1+1}|m]+\bE\left[[\alpha\frac{m+1}{n_1+1}-\frac{n_0}{n_1+1}]_+|m\right],\notag\\
&=\bE[\frac{1}{n_1+1}|m]+\bE\left[[\alpha\frac{m+1-n_0}{n_1+1}-(1-\alpha)\frac{n_0}{n_1+1}]_+|m\right],\notag\\
&=\bE[\frac{1}{n_1+1}|m]+(1-\alpha)\bE\left[[\frac{\alpha}{1-\alpha}-\frac{n_0}{n_1+1}]_+|m\right].
\end{align}
Notice that conditional on $m$, $X_i$ falls at $0$ or $X_{n+1}$ following an independent Bernoulli law: 
\[
X_i = \left\{\begin{array}{ll} 0 &\mbox{ w.p. } \frac{q_0}{q_0+q_1} = \alpha,\\
X_{n+1}&\mbox{ w.p. } \frac{q_1}{q_0+q_1} =1- \alpha.
\end{array}\right.
\]
 From direct calculations, we obtain that
\begin{align}
\label{eq:counter2_eq2}
\bE[\frac{1}{n_1+1}|m] &= \sum_{n_1=0}^m\frac{1}{n_1+1}\frac{m!}{n_1!(m-n_1)!}(1-\alpha)^{n_1}\alpha^{m-n_1}\notag\\
&=\frac{1}{(m+1)(1-\alpha)}\sum_{n_1=1}^{m+1} \frac{(m+1)!}{n_1!(m+1-n_1)!}(1-\alpha)^{n_1}\alpha^{m+1-n_1}\leq \frac{1}{m(1-\alpha)}.
\end{align}
Also, we have
\begin{align}
\label{eq:counter2_eq3}
&\bE\left[[\frac{\alpha}{1-\alpha}-\frac{n_0}{n_1+1}]_+|m\right]\notag\\
\overset{(c)}{=}&\frac{\alpha}{1-\alpha}\bP(n_0\leq \alpha(m+1)|m)-\sum_{n_0=1}^{n_0\leq \alpha(m+1)}\frac{n_0}{m-n_0+1}\frac{m!}{n_0!(m-n_0)!}\alpha^{n_0}(1-\alpha)^{m-n_0}\notag\\
=&\frac{\alpha}{1-\alpha}\left(\bP(n_0\leq \alpha(m+1)|m)-\sum_{n_0=0}^{n_0\leq \alpha(m+1)-1}\frac{m!}{n_0!(m-n_0)!}\alpha^{n_0}(1-\alpha)^{m-n_0}\right)\notag\\
=&\frac{\alpha}{1-\alpha}\left(\bP(n_0\leq \alpha(m+1)|m)-\bP(n_0\leq \alpha(m+1)-1|m)\right)\notag\\
=&\frac{\alpha}{1-\alpha}\bP(n_0 = \underbrace{\lfloor\alpha(m+1)\rfloor}_{n_*} |m)=\frac{\alpha}{1-\alpha}{m \choose n_*} \alpha^{n_*}(1-\alpha)^{m-n_*}.
\end{align}
We  now use the Stirling's approximation:
\[
 \sqrt{2\pi k}(\frac{k}{e})^ke^{\frac{1}{12k+1}}\leq k!\leq \sqrt{2\pi k}(\frac{k}{e})^ke^{\frac{1}{12k}},\; \mbox{ for all }k\geq 1.
\]
Plug the Stirling's approximation into (\ref{eq:counter2_eq3}), there exist a constant $C>0$ such that when $m\geq C$, we have:
\begin{align}
\label{eq:counter2_eq4}
\bE\left[[\frac{\alpha}{1-\alpha}-\frac{n_0}{n_1+1}]_+|m\right]&\leq \frac{\alpha}{1-\alpha}\exp(\frac{1}{12m})\sqrt{\frac{m}{2\pi (n^*)(m-n^*)}}(\alpha\frac{m}{n^*})^{n^*}\left((1-\alpha)\frac{m}{m-n^*}\right)^{m-n^*}\notag\\
& \overset{(c)}{\leq}  \frac{\alpha}{1-\alpha}\exp(\frac{1}{12m})\sqrt{\frac{m}{2\pi (m\alpha-1)(m(1-\alpha)-1)}}(\frac{m\alpha}{m\alpha-1})^{n^*}\left(\frac{(1-\alpha)m}{m(1-\alpha)-1}\right)^{m-n^*}\notag\\
&\leq C\sqrt{\frac{1}{m}}(1+\frac{2}{m})^m\leq \frac{Ce^2}{\sqrt{m}},
\end{align}
where we have used the fact that $m\alpha+1\leq n^* \geq \alpha m-1$ at step (c). Notice that $m$ itself follows a binomial distribution with $n$ trials and successful rate $(q_1+q_0)$. Apply the Chernoff bound, we have
\begin{equation}
\label{eq:counter2_eq5}
\bP(m\leq \frac{(q_1+q_0)n}{2})\leq \exp(-\frac{n(q_1+q_0)}{8}).
\end{equation}
For any constant $p\geq 1$,  $n(q_1+q_0)\rightarrow\infty$. Combine it with (\ref{eq:counter2_eq1}), (\ref{eq:counter2_eq2}), (\ref{eq:counter2_eq4}) and (\ref{eq:counter2_eq5}),  there exist a constant $C>0$, such that for all $X_{n+1}\in \{\pm e_j, j = 1,\ldots, p\}$, we have
\begin{align*}
&\bP(V_{n+1}\leq Q(\alpha; \hat\mF)|X_{n+1})\\
\leq &\bP(\{V_{n+1}\leq Q(\alpha; \hat\mF)|X_{n+1}\}\cap \{m\geq \frac{(q_1+q_0)n}{2}\})+ \bP(m\geq \frac{(q_1+q_0)n}{2}|X_{n+1})\leq C\sqrt{\frac{1}{(q_1+q_0)n}}.
\end{align*}
Marginalize over $X_{n+1}$, we reach the desired result: there exists a sufficiently large constant $C$, such that
\[
\bP(V_{n+1}\leq Q(\alpha; \hat\mF))\leq \bP(X_{n+1})\frac{C}{\sqrt{(q_0+q_1)n}}+\bP(X_{n+1}=0)\leq \frac{C}{\sqrt{(q_0+q_1)n}}+q_0\rightarrow q_0.
\]
\end{proof}
\subsection{Proof of Theorem \ref{thm:general1}}
\begin{proof}
Define 
\[
\mathcal{T} \coloneqq \left\{\{Z_i, i = 1,\ldots, n+1\}= \{z_i\coloneqq (x_i, y_i), i = 1,\ldots, n+1\}\right\}.
\]
Let $\sigma$ be a permutation of numbers $1, 2,\ldots, n+1$ that specifies how the values are assigned, e.g., $Z_i$ takes value $z_{\sigma_i}$.  Since $V(.;\mZ)$ and $H(.,.;\mX)$ are fixed conditional on $\mathcal{T}$, we can set $v_{\sigma_i}^* = Q(\tilde\alpha; \sum_{j=1}^n p^H_{\sigma_i j}\delta_{v_{j}})$ as the realized empirical quantile at $\tilde\alpha$ for $\hat\mF_i$ given a particular permutation ordering $\sigma$.  Hence, for any given $\tilde\alpha\in \Gamma$, conditional $\mathcal{T}$ and the permutation ordering $\sigma$, we have
 \begin{align}
 \label{eq:proof_thm1eq1}
\sum^{n+1}_{i=1}\mathbbm{1}_{V_i \leq Q(\tilde\alpha;\hat\mF_i)}|\mathcal{T},\sigma& = \sum^{n+1}_{i=1}\mathbbm{1}_{v_{\sigma_i} \leq v^*_{\sigma_i}}=\sum^{n+1}_{i=1}\mathbbm{1}_{v_i \leq v^*_{i}}.
\end{align}
In other words, the achieved value for the left side of (\ref{eq:proof_thm1eq1}) or  Theorem \ref{thm:general1} (\ref{eq:goal1general}) remains the same for all $\sigma$.  Since $\Gamma$  is fixed conditional on $\mathcal{T}$, the smallest value in $\Gamma$ satisfying (\ref{eq:goal1general}) is also fixed conditional on $\mathcal{T}$, by Lemma  \ref{lem:conditional}, we obtain that
\[
\bP\{V_{n+1} \leq Q(\tilde\alpha; \hat \mF_{n+1}) |\mathcal{T}\}= \bE\left\{\frac{1}{n+1}\sum^{n+1}_{i=1}\mathbbm{1}_{v_i \leq v^*_i}|\mathcal{T}\right\}=\frac{1}{n+1}\sum^{n+1}_{i=1}\mathbbm{1}_{v_i \leq v^*_i}\geq\alpha
\]
Marginalize over $\mathcal{T}$, we have
\begin{equation}
\label{eq:result1general}
\bP\{V_{n+1}\leq  Q(\tilde\alpha; \hat \mF_{n+1}) \} \geq \alpha.
\end{equation}
By Lemma \ref{lem:infequiv}, equivalently, we also have
\begin{equation}
\label{eq:result2general}
\bP\{V_{n+1}\leq Q(\tilde\alpha; \hat \mF)\} \geq \alpha.
\end{equation}
\end{proof}

\subsection{Proof of Theorem \ref{thm:general2}}
Define
\[
\mathcal{T} \coloneqq \left\{\{Z_i, i = 1,\ldots, n+1\}= \{z_i\coloneqq (x_i, y_i), i = 1,\ldots, n+1\}\right\}.
\]
By (\ref{eq:proof_thm1eq1}) and the fact that $\Gamma$ is fixed conditional on $\mathcal{T}$, we know that $\tilde\alpha_1$, $\tilde\alpha_2$ and $\alpha_1$, $\alpha_2$ are fixed conditional on $\mathcal{T}$.  As a result,  when $\tilde\alpha = \left\{\begin{array}{cc}\tilde \alpha_1&w.p.\;\frac{\alpha - \alpha_2}{\alpha_1-\alpha_2}\\
\tilde\alpha_2 & w.p.\;\frac{\alpha_1 - \alpha}{\alpha_1-\alpha_2}
\end{array}\right.$, and it is  independent of the data conditional on $\mathcal{T}$. Apply Lemma \ref{lem:conditional}, we have
\begin{align*}
\bP\{V_{n+1} \leq Q(\tilde\alpha; \hat \mF_{n+1}) |\mathcal{T}\}&= \bE\left\{\frac{1}{n+1}\sum^{n+1}_{i=1}\mathbbm{1}_{v_i \leq v^*_i}|\mathcal{T}\right\} \\
& = \alpha_1\frac{\alpha-\alpha_2}{\alpha_1 -\alpha_2}+\alpha_2\frac{\alpha_1 - \alpha}{\alpha_1 - \alpha_2} = \alpha.
\end{align*}
Marginalizing over $\mathcal{T}$, we have
\[
\bP\{V_{n+1} \leq Q(\tilde\alpha; \hat \mF_{n+1})\}= \alpha.
\]
By Lemma \ref{lem:infequiv}, equivalently, we have
\[
\bP\{V_{n+1} \leq Q(\tilde\alpha; \hat \mF)\}= \alpha.
\]

\subsection{Proof of Lemma \ref{lem:split1}}
\begin{proof}
As a direct application of Theorem \ref{thm:general1} and Crorllary \ref{cor:generalC}, we obtain that
\[
\bP\left\{V_{n+1} \in C_V(X_{n+1}) \right\}\geq \alpha,\; \bP\left\{Y_{n+1} \in C(X_{n+1}) \right\}\geq \alpha.
\]
The fact that  $C_V(X_{n+1})$ is an interval comes directly from Lemma \ref{lem:monotone} (iii).
\end{proof}
\subsection{Proof of Lemma \ref{lem:alg1} }
\begin{proof}
\;
\begin{itemize}
\item Proof of part 1: By definition, $V_{n+1} = v\in C_V(X_{n+1})$ iff (if and only if) the smallest value $\tilde\alpha\in \Gamma$ that makes (\ref{eq:goal1}) hold is greater than $\sum_{V_i < v} p^H_{n+1,i}\in \Gamma$. That is, $v\in C_V(X_{n+1})$ iff 
\begin{equation}
\label{eq:proof_alg1}
\frac{1}{n+1}\sum_{i=1}^n \mathbbm{1}_{V_i\leq Q(\sum_{V_i < v} p^H_{n+1,i}; \hat\mF_i(v))}< \alpha.
\end{equation}
\begin{itemize}
\item[(a)]  When $v= \overline{V}_k$ for some $1\leq k\leq n+1$,  $\sum_{V_i < \overline{V}_k} p^H_{n+1,i} = \tilde\theta_k$ by definition. Hence $v\in C_V(X_{n+1})$ iff
\begin{equation}
\label{eq:proof_alg1}
\frac{1}{n+1}\sum_{i=1}^n \mathbbm{1}_{V_i \leq Q(\tilde\theta_k; \hat\mF_i(\overline{V}_{k}))} < \alpha.
\end{equation}
\item[(b)] When $v\in (\overline{V}_{\ell(k)}, \overline{V}_{k})$ for some $1\leq k\leq n+1$,  $\sum_{V_i < v} p^H_{n+1,i} = \sum_{V_i < \overline{V}_{k}} p^H_{n+1,i}=\tilde\theta_k$. Hence $v\in C_V(X_{n+1})$ iff
\begin{equation}
\label{eq:proof_alg2}
\frac{1}{n+1}\sum_{i=1}^n \mathbbm{1}_{V_i \leq Q(\tilde\theta_k; \hat\mF_i(v))} < \alpha.
\end{equation}
A key observation is that the status of event $ \{V_i \leq Q(\tilde\theta_k; \hat\mF_i(v))\}$ does not change as we vary $v\in [\overline{V}_{\ell(k)}, \overline{V}_{k})$. That is, for all $1\leq i\leq n$, we have 
\[
J_{ik}(v)\coloneqq \{V_i \leq Q(\tilde\theta_k; \hat\mF_i(v))\}= \{V_i \leq Q(\tilde\theta_k; \hat\mF_i(\overline{V}_{\ell(k)})\}\coloneqq J_{ik}.
\]
This can be easily verified:
\begin{itemize}
\item  [$\diamond$] If $V_i < \overline{V}_k$, we have $V_i \leq \overline{V}_{\ell(k)} <v$, and
\[
\{V_i \leq Q(\tilde\theta_k; \hat\mF_i(v))\}=\{\tilde\theta_k > \theta_i\} =  \{V_i \leq Q(\tilde\theta_k; \hat\mF_i(\overline{V}_{\ell(k)}))\}
\]
\item [$\diamond$] If $V_i \geq  \overline{V}_k$, then $V_i > v > \overline{V}_{\ell(k)}$, and we have
\[
\{V_i \leq Q(\tilde\theta_k; \hat\mF_i(v))\}=\{\tilde\theta_k > \theta_i+p^H_{i,n+1}\} =  \{V_i \leq Q(\tilde\theta_k; \hat\mF_i(\overline{V}_{\ell(k)}))\},
\] 
\end{itemize}
Hence, we obtain
\begin{equation}
\label{eq:proof_alg3}
\frac{1}{n+1}\sum_{i=1}^n \mathbbm{1}_{V_i \leq Q(\tilde\theta_k; \hat\mF_i(\overline{V}_{l(k)}))} < \alpha.
\end{equation}
\end{itemize}
Combine part (a) and part (b), and the fact that $\overline{V}_{\ell(k)}\leq \overline{V}_{k}$ and $Q(\tilde\theta_k; \hat\mF_i(v))$ is non-decreasing in $v$ (Lemma \ref{lem:monotone}), we immediately reach the desired result that
\[
\bar{C}^V(X_{n+1}) = \{v: v \leq  Q(\tilde\theta_{k^*}; \hat\mF)\},
\]
where $k^*$ is the largest value of $k$ such that (\ref{eq:proof_alg3}) holds.
\item Proof of  part 2:  As we increase $k$, both $\overline{V}_{l(k)}$ and $\tilde \theta_{k}$ are non-decreasing,  hence, $Q(\tilde\theta_{k}; \hat{\mF}_i(\overline{V}_{\ell(k)}))$ is non-decreasing in $k$. Thus, $J_{ik} = \{V_i\leq Q(\tilde\theta_k; \hat\mF_i(\bar{V}_{\ell(k)}))\}$ is a monotone event in $k$: for all $k'\geq k$,  we have $J_{ik}\subseteq J_{ik'}$.  Consequently, suppose $k^*_i$ is when $J_{ik}$ first holds, then $\mathbbm{1}_{J_{ik}} = 1$ iff $k\geq k^*_i$.  We can divide $J_{ik}$ into two subsets:
\begin{align}
\label{eq:proof_alg4}
J_{ik} &= \left(\{V_i > \bar{V}_{\ell(k)}\}\cap \{V_i\leq Q(\tilde\theta_k; \hat\mF_i(\bar{V}_{\ell(k)}))\}\right)\cup \left(\{V_i \leq \bar{V}_{\ell(k)}\}\cap \{V_i\leq Q(\tilde\theta_k; \hat\mF_i(\bar{V}_{\ell(k)}))\}\right)\notag\\
& \overset{(a)}{=}  \underbrace{\left(\{\ell(i) \geq \ell(k)\}\cap \{\theta_i+p^H_{i,n+1} <\tilde\theta_k\}\right)}_{J^1_{ik}}\cup   \underbrace{\left(\{\ell(i) < \ell(k)\}\cap \{\theta_i <\tilde\theta_k\}\right)}_{J^2_{ik}}.
\end{align}
At step (a), we have used the fact that
\[
V_i > \overline{V}_{\ell(k)}\Leftrightarrow \ell(i)\geq \ell(k),
\]
and that 
\begin{itemize}
\item  when  $V_i > \overline{V}_{\ell(k)}$, we have $\sum_{1\leq j\leq n+1:V_j < V_i}p_{ij}^H = \theta_i + p^H_{i,n+1}$ when $V_{n+1} = \bar{V}_{\ell(k)}$ Hence, 
\[
V_i\leq Q(\tilde\theta_k; \hat\mF_i(\bar{V}_{l(k)}))\Leftrightarrow \theta_i+p^H_{i,n+1} <\tilde\theta_k.
\]
\item  when  $V_i \leq  \overline{V}_{\ell(k)}$,  we have $\sum_{1\leq j\leq n+1:V_j < V_i}p_{ij}^H = \theta_i$ when $V_{n+1} = \bar{V}_{\ell(k)}$.  Hence, 
\[
V_i\leq Q(\tilde\theta_k; \hat\mF_i(\bar{V}_{l(k)}))\Leftrightarrow \theta_i <\tilde\theta_k.
\]
\end{itemize}
We now consider when $J_{ik}$ turns true for samples from categories $A_1$, $A_2$ and $A_3$.
\begin{itemize}
\item For $i\in A_1$, by the definition of $A_1$ and (\ref{eq:proof_alg4}), we know that $J_{ik}$ is true at $k = i$ and $k_i^*\leq i$, and $J_{ik} =J_{ik}^1= \{\theta_i+p^H_{i,n+1} < \tilde\theta_{k}\}$.
\item For $i\in A_2\cup A_3$, since $i\notin A_1$ and $k_i^* > i$, $J_{ik}^1$ fails to hold for all $k$. Hence, $J_{ik}$ holds when $J_{ik}^2$ holds. 
\begin{itemize}
\item[$\diamond$] When $i\in A_2$: since $\theta_i \geq \tilde\theta_i$, in order for $\theta_i < \tilde\theta_k$ to hold,   by definition, we must have $ \sum_{j\leq l(k)}p^H_{n+1,j}=\tilde\theta_k>\tilde\theta_i = \sum_{j\leq l(i)}p^H_{n+1,j}$, which automatically guarantees that $l(k) > l(i)$. As a result,  for $i\in A_2$, we have $J_{ik} = J_{ik}^2 = \{\theta_i < \tilde\theta_k\}$.
\item[$\diamond$] When $i\in A_3$: in order to have $l(k)>l(i)$, we automatically $\tilde\theta_k \geq \tilde\theta_i>\theta_i$ for samples in $A_3$. Thus,   for $i\in A_3$, we have $J_{ik} = J_{ik}^2 = \{l(i)<l(k)\}$. 
\end{itemize}
\end{itemize}
Combine them together, we have
\begin{align*}
S(k)&=\sum_{i=1}^{n}\frac{1}{n+1}\mathbbm{1}_{V_i \leq Q(\tilde\theta_{k}; \hat{\mathcal{F}}_i(\overline{V}_{\ell(k)}))}  \\
&=\frac{1}{n+1}\left(\sum_{i\in A_1}\mathbbm{1}_{J_{ik}^1}+ \sum_{i\in A_2}\mathbbm{1}_{J_{ik}^2}+ \sum_{i\in A_3}\mathbbm{1}_{J_{ik}^2}\right)\\
&=\frac{1}{n+1}\left(\sum_{i\in A_1}\mathbbm{1}_{\{\theta_i+p^H_{i,n+1} < \tilde\theta_{k}\}}+ \sum_{i\in A_2}\mathbbm{1}_{\theta_i < \tilde\theta_k}+ \sum_{i\in A_3}\mathbbm{1}_{l(i)<l(k)}\right).
\end{align*}
We have proved the second part of Lemma \ref{lem:alg1}.
\end{itemize}
\end{proof}
\section*{Local coverage properties of LCP}
\subsection{Proof of Theorem \ref{thm:AsymptoticConditional}}

\begin{proof} We first prove the convergence from $\tilde\alpha(v)$ to $\alpha$ in (\ref{eq:thmAsy3}) and then show that the achieved coverage levels converge to the nominal level for both $\tilde\alpha = \alpha$ and $\tilde\alpha = \tilde\alpha(v)$ as described in Lemma \ref{lem:split1}. Define  $I_{i}= \frac{\sum_{j=1, j\neq i}^n P_{V|X_j}(V_i) H_{ij}}{B_i}$ and $R_i = \sum_{j=1, j\neq i}^n\frac{H_{ij}(\mathbbm{1}_{V_j < V_i} - P_{V|X_j}(V_i))}{B_i}$ for all $i=1,\ldots, n+1$.
\begin{enumerate}
\item {\textbf{Proof of (\ref{eq:thmAsy3}):}}
For $i=1,\ldots, n$, define $B_i =  \sum_{j=1}^{n+1} H_{ij}$ and, for any $\tilde\alpha \in [0, 1]$ and $v\in \real$,  define
\[
J_i(v, \tilde\alpha) \coloneqq \{V_i \leq Q(\tilde\alpha; \hat\mF(v))\}=\{\tilde\alpha > \frac{\sum_{j\leq n:V_j < V_i} H_{ij}+\mathbbm{1}_{v<V_i}}{B_i}\}.
\] 
$J_i(v, \tilde\alpha)$ is the event for wether sample $i$ contributes to the left side of Lemma \ref{lem:split1} (\ref{eq:goal1}). We can define a subset event $\underline{J}_{i}(\tilde\alpha)\subseteq J_i(v, \tilde\alpha)$ for all $v$ values for all $v$. Decompose the condition of $J_{i}(v, \tilde\alpha)$ as below:
\begin{align}
\label{eq:proof_asym_upper1}
\frac{\sum_{j\leq n:V_j < V_i} H_{ij}+\mathbbm{1}_{v<V_i}}{B_i}&\leq \frac{\sum_{j\leq n:V_j < V_i} H_{ij}}{B_i}+\frac{1}{B_i}=I_i + R_i+\frac{1}{B_i}.
\end{align}
Set $G = \left\{i\in \{1,\ldots, n\}: B_i\geq \frac{1}{2eL}n h_n^\beta, |R_i|\leq \sqrt{\frac{2eL\ln n}{n h_n^\beta}}\right\}$. By Lemma \ref{lem:alpha_bound} (i), there exists a constant $C>0$, such that for all $i \in G$:
{\small
\begin{align}
\label{eq:proof_asym_bound1}
\frac{B_i-1-H_{i,n+1}}{B_i}\in [1-\frac{4eL}{nh_n^\beta},1],\quad |I_i - \frac{B_i-1-H_{i,n+1}}{B_i}P_{V|X_i}(V_i)|\leq C h_n \ln (h_n^{-1}).
\end{align}
}
Combine (\ref{eq:proof_asym_upper1}) with  (\ref{eq:proof_asym_bound1}), there exist a constant $C>0$ such that for all $i \in G$, we have
\begin{equation}
\label{eq:proof_asym_upper2}
\underline{J}_i(\tilde\alpha)\coloneqq \left\{\tilde\alpha > P_{V|X_i}(V_i)+C \left(h_n \ln (h_n^{-1})+\sqrt{\frac{\ln n}{n h_n^\beta}}\right)\right\}\subseteq J_i(v, \tilde\alpha),\;\mbox{for all } v.
\end{equation}
We can also define a superset event $\bar{J}_i(\tilde\alpha)\supseteq J_i(v, \tilde\alpha)$ for all $v$ values:
\begin{align}
\label{eq:proof_asym_lower1}
\frac{\sum_{j\leq n:V_j < V_i} H_{ij}+\mathbbm{1}_{v<V_i}}{B_i}&\geq \frac{\sum_{j\leq n:V_j < V_i} H_{ij}}{B_i}=I_i + R_i.
\end{align}
Combine (\ref{eq:proof_asym_lower1}) with  (\ref{eq:proof_asym_bound1}), there exists a constant $C>0$ such that for all $i \in G$, we have
\begin{equation}
\label{eq:proof_asym_lower2}
J_i(v, \tilde\alpha)\subseteq \bar{J}_i(\tilde\alpha)\coloneqq \left\{\tilde\alpha > P_{V|X_i}(V_i)-C \left(h_n \ln (h_n^{-1})+\sqrt{\frac{\ln n}{n h_n^\beta}}\right)\right\} ,\;\mbox{for all } v.
\end{equation}
Hence, we can then upper and lower bound the left side of (\ref{eq:goal1}) using $\bar{J}_i(\tilde\alpha)$ and $\underline{J}_i(\tilde\alpha)$:
\begin{align}
&\frac{1}{n+1}\sum_{i=1}^{n+1} J_i(v, \tilde\alpha) \leq \frac{1}{n+1}+\frac{1}{n+1}\sum_{i\in G} \bar{J}_i(\tilde\alpha)+\frac{|G^c|}{n+1} \label{eq:proof_asym_upper3},\\
&\frac{1}{n+1}\sum_{i=1}^{n+1} J_i(v, \tilde\alpha) \geq\frac{1}{n+1}\sum_{i\in G} \underline{J}_i(\tilde\alpha) \label{eq:proof_asym_lower3}.
\end{align}
Set $W_i = P_{V|X_i}(V_i)$, which is i.i.d generated from $\operatorname{Unif}[0,1]$  when $V|X_i$ is a continuous variable.  By Lemma \ref{lem:alpha_bound}, we know that
\begin{align}
\label{eq:proof_G_bound}
\bP\{|G^c|=0\}&\leq \bP\{\min_{i=1}^n B_i \leq \frac{n h_n^\beta}{2eL}\}+\bP\{\exists i\in \{1,\ldots, n\}: B_i>\frac{n h_n^\beta}{2eL}, |R_i|\leq \sqrt{\frac{2eL\ln n}{ n h_n^\beta}}\}\notag\\
& \leq  \bP\{\min_{i=1}^n B_i \leq \frac{n h_n^\beta}{2eL}\}+\bP\{\max_{1\leq i\leq n}|R_i|\leq \sqrt{\frac{\ln n}{B_i}}\}\notag\\
& \leq n\exp(-\frac{nh_n^\beta}{8L})+n\times \frac{1}{n^2}\rightarrow 0.
\end{align}
When $\{|G^c|=0\}$ holds:
\begin{itemize}
\item When $\tilde\alpha$ makes (\ref{eq:goal1}) hold, by (\ref{eq:proof_asym_upper3}),  we must have
\begin{align}
\label{eq:proof_alpha_low1}
&\frac{1}{n+1}\left(1+ \sum_{i=1}^n\bar{J}_i(\tilde\alpha)\right) \geq \alpha\Rightarrow \tilde\alpha \geq Q(\frac{n+1}{n}\alpha-\frac{1}{n}; \frac{1}{n}\sum_{i=1}^n\delta_{W_i})-C \left(h_n \ln (h_n^{-1})+\sqrt{\frac{\ln n}{n h_n^\beta}}\right).
\end{align}
\item   By (\ref{eq:proof_asym_lower3}), $\tilde\alpha$ makes (\ref{eq:goal1}) hold as long as
\begin{align*}
&\frac{1}{n+1}\sum_{i=1}^n \underline{J}_i(\tilde\alpha)\geq \alpha\Rightarrow \tilde\alpha \geq Q(\frac{n+1}{n}\alpha; \frac{1}{n}\sum_{i=1}^n\delta_{W_i})+C \left(h_n \ln (h_n^{-1})+\sqrt{\frac{\ln n}{n h_n^\beta}}\right).
\end{align*}
Further, since $\Gamma$ includes all possible empirical CDF values from weighted distribution $\hat\mF_i$ for $i=1,\ldots, n+1$ under all possible ordering of $V_1,\ldots, V_{n+1}$. Let $B_{\max} = \max_{i=1}^{n+1} B_i$. The differences between two adjacent values in $\Gamma$ is upper bounded by $\frac{1}{B_{\max}}\leq \frac{2eL}{nh_n^\beta}$. Hence, there exists a constant $C>0$ such that the smallest value in $\Gamma$ that makes (\ref{eq:goal1}) is upper bounded by
\begin{align}
\label{eq:proof_alpha_up1}
\tilde\alpha\leq Q(\frac{n+1}{n}\alpha; \frac{1}{n}\sum_{i=1}^n\delta_{W_i})+C \left(h_n \ln (h_n^{-1})+\sqrt{\frac{\ln n}{n h_n^\beta}}\right).
\end{align}
\end{itemize}
The bounds (\ref{eq:proof_alpha_low1}) and (\ref{eq:proof_alpha_up1}) hold for all $V_{n+1}=v$. By Dvoretzky–Kiefer–Wolfowitz inequality and the fact that $W_i\sim\operatorname{Unif}[0,1]$, there exists a constant $C>0$ such that
\begin{align}
\label{eq:proof_asy_bound2}
\bP(\max_t|Q(t; \sum_{i=1}^n\delta_{W_i})-t|\leq \sqrt{\frac{\ln n}{n}})\leq \frac{C}{n^2}.
\end{align}
Combine \ref{eq:proof_alpha_low1},  (\ref{eq:proof_asy_bound2}) and (\ref{eq:proof_G_bound}), there exist a constant $C>0$, such that
{\small
\begin{align}
\label{eq:proof_alpha_low2}
\bP\left\{|\min_{v_{n+1}}\tilde\alpha(v_{n+1})-\alpha|< C \left(h_n \ln (h_n^{-1})+\sqrt{\frac{\ln n}{n h_n^\beta}}\right)\right\}\geq \frac{C}{n^2}+\bP(|G^c|>0)\rightarrow 0.
\end{align}
}
Since $C \left(h_n \ln (h_n^{-1})+\sqrt{\frac{\ln n}{n h_n^\beta}}\right)\rightarrow 0$, this concludes our proof.
\item {\textbf{Proofs of (\ref{eq:thmAsy1}) and (\ref{eq:thmAsy2}):}} By definition, for any given $\tilde\alpha$, $V_{n+1}\leq Q(\tilde\alpha; \hat\mF)$ if and only if
\begin{equation}
\label{eq:proof_alpha_coverage1}
\sum_{i=1}^n \frac{H(X_{n+1},X_i)}{B_{n+1}}\mathbbm{1}_{V_i < V_{n+1}}=I_{n+1}+R_{n+1}< \tilde\alpha.
\end{equation}
Define $G= \{B_{n+1}\geq \frac{n h_n^\beta}{2eL}, |R_{n+1}|\leq \sqrt{\frac{2eL\ln n}{ n h_n^\beta}}\}$. When $G$ holds, following the same routine as bounding $J(\tilde\alpha, v)$ with $\bar{J}(\tilde\alpha)$ and $\underline{J}(\tilde\alpha)$, we can lower and upper bound the left side of (\ref{eq:proof_alpha_coverage1}) using Lemma \ref{lem:alpha_bound} (i): there exists a constant $C>0$, such that
\begin{align}
&I_{n+1}+R_{n+1}\leq P_{V|X_{n+1}}(V_{n+1})+C \left(h_n \ln (h_n^{-1})+\sqrt{\frac{\ln n}{n h_n^\beta}}\right).\\
&I_{n+1}+R_{n+1}\geq  P_{V|X_{n+1}}(V_{n+1})-C \left(h_n \ln (h_n^{-1})+\sqrt{\frac{\ln n}{n h_n^\beta}}\right).
\end{align}
$W_{n+1}=P_{V|X_{n+1}}(V_{n+1})\sim \operatorname{Unif}[0,1]$ since $V|X_{n+1}$  is a continuous variable.  By Lemma \ref{lem:alpha_bound} (i) and (ii), $\bP(G^c)\rightarrow 0$. Hence, for any given $\tilde\alpha$, there exists a constant $C>0$, such that
\begin{align}
& \bP(I_{n+1}+R_{n+1} < \tilde\alpha)\leq \tilde\alpha+C\left(h_n\ln (h_n^{-1})+\frac{1}{(nh_n^{\beta})^{1\slash 3}}\right)+\bP\left\{G^c\right\}\rightarrow \tilde\alpha,\\
& \bP(I_{n+1}+R_{n+1} < \tilde\alpha)\geq \tilde\alpha-C\left(h_n\ln (h_n^{-1})+\frac{1}{(nh_n^{\beta})^{1\slash 3}}\right)\rightarrow \tilde\alpha.
\end{align}
Consequently, when $\tilde\alpha = \alpha$ or $\tilde\alpha = \tilde\alpha(v)\rightarrow \alpha$ for all $v$ in probability as described in (\ref{eq:thmAsy3}), we achieve an asymptotic conditional coverage at level $\alpha$.
\end{enumerate}
\end{proof}

\subsection{Proof of Theorem \ref{thm:ApproximateConditional}}
\begin{proof}
We use the result from \cite{barber2019conformal}  which extends CP to the setting with covariate shift:
\begin{proposition}[B.1 (\citet{barber2019conformal}, Corollary 1)]
\label{prop:shift}
For any fixed $x_0$. Set $w_{x_0}(.) = \frac{d \tilde \mP^{x_0}_X}{d \mP_X}$ and  $p^{x_0}_i(x)=\frac{w_{x_0}(X_i)}{\sum_{j=1}^n w_{x_0}(X_i)+w_{x_0}(x)}$ for $i=1,\ldots, n$, and $p^{x_0}_{n+1}(x)=\frac{w_{x_0}(x)}{\sum_{j=1}^n w_{x_0}(X_i)+w_{x_0}(x)}$. Then,
\[
\bP(V(X_{n+1}, Y_{n+1}) \leq Q(\alpha; \sum_{i=1}^{n+1} p^{x_0}_i(X_{n+1})\delta_{\overline{V}_i} \}) \geq\alpha.
\]
\end{proposition}
In our setting, $w_{x_0}(x) \propto H(x_0, x)$.  As a direct application of Proposition \ref{prop:shift}, when $(\tilde X,\tilde Y)$ is distributed from $\tilde\mP^{X_{n+1}}_{XY}$, we have 
\[
\bP\{V(\tilde X, \tilde Y)\leq Q(\alpha; \sum_{i=1}^{n+1} p^{x_0}_i(\tilde X)\delta_{\overline{V}_i})|X_{n+1}=x_0\}\geq \alpha.
\]
Since the $H(x_0, x_0)\geq H( x_0, \tilde X)$ by definition, the distribution $\hat\mF$ dominates the distribution  $\sum_{i=1}^{n+1} p^{x_0}_i(\tilde X)\delta_{\overline{V}_i}$: given $X_{n+1}=x_0$, for any $\alpha$, we have
\begin{align*}
Q(\alpha;\hat\mF) &=Q(\alpha;\sum_{i=1}^{n+1}\frac{H(x_0,X_i)}{\sum_{j=1}^{n+1}H(x_0,X_j)}\delta_{\bar{V}_i})\\
&\geq Q(\alpha;\sum_{i=1}^{n}\frac{H(x_0,X_i)}{\sum_{j=1}^{n}H(x_0,X_j)+H(x_0, \tilde X)}\delta_{\bar{V}_i}+\frac{H(x_0,\tilde X)}{\sum_{j=1}^{n}H(x_0,X_j)+H(x_0, \tilde X)}\delta_{\bar{V}_i})= Q(\alpha; \sum_{i=1}^{n+1} p^{x_0}_i(\tilde X)\delta_{\overline{V}_i}).
\end{align*}
Hence, we have
\[
\bP\{V(\tilde X, \tilde Y)\leq  Q(\alpha;\hat\mF)|X_{n+1}=x_0\}\geq \alpha,\;\mbox{for all } x_0.
\]
Next, we turn to the achieved coverage using $\tilde C(X_{n+1})$. By construction, we have
\begin{align*}
\{\tilde Y \in \tilde C(X_{n+1})\}&=\{V(X_{n+1}, \tilde Y)\leq Q(\alpha; \hat\mF)+\varepsilon(X_{n+1})\}\\
&\supseteq \{V(\tilde X, \tilde Y)\leq Q(\alpha; \hat\mF)\}.
\end{align*}
Consequently, we obtain
\[
\bP(\tilde Y \in \tilde C(X_{n+1})|X_{n+1}=x_0) \geq \bP(V(\tilde X) \leq Q(\alpha;\hat\mF))|X_{n+1}=x_0) \geq \alpha.
\]
\end{proof}

\section{Proof of Lemmas in the Appendix}
\label{app:proofapp}
\subsection{Proof of Lemma \ref{lem:infequiv}}
\begin{proof}
By definition, we know 
\[
V_{n+1} \leq Q(\alpha; \sum^n_{i=1} w_i \delta_{V_i} + w_{n+1}\delta_{V_{n+1}}) \Rightarrow V_{n+1} \leq Q(\alpha; \sum^n_{i=1} w_i \delta_{V_i} + w_{n+1}\delta_\infty).
\]
To show that Lemma \ref{lem:infequiv} holds, we only need to show that,
\[
V_{n+1} >Q(\alpha; \sum^n_{i=1} w_i \delta_{V_i} + w_{n+1}\delta_{V_{n+1}})\Rightarrow V_{n+1} > Q(\alpha; \sum^n_{i=1} w_i \delta_{V_i} + w_{n+1}\delta_\infty).
\]
Let $Q(\alpha;\sum_{i=1}^n w_i\delta_{V_i}+w_{n+1}V_{n+1}) = V_{i^*}$ for some index $1\leq i^*\leq n+1$. When $V_{n+1} >V_{i^*}$, we must have $V_{i^*}<\infty$. By definition:
 \begin{align*}
 & \alpha\geq \sum_{i=1}^{n+1}w_i\mathbbm{1}_{V_i\leq V_{i^*}} = \sum_{i=1}^nw_i\mathbbm{1}_{V_i\leq V_{i^*}} =\sum_{i=1}^nw_i\mathbbm{1}_{V_i\leq V_{i^*}} +w_{n+1}\mathbbm{1}_{\infty\leq V_{i^*}}\\\Rightarrow &Q(\alpha;\sum_{i=1}^n w_i \mathbbm{1}_{V_i}+w_{n+1}\delta_{\infty} )\leq V_{i^*}< V_{n+1}.
 \end{align*}
\end{proof}

\subsection{Proof of Lemma \ref{lem:monotone}}
\begin{proof} 
We can prove Lemma  \ref{lem:monotone} with elementary calculus  arguments.
\begin{itemize}
\item[(i)]  Given $V_{1}, \ldots, V_{n+1}$, $Q(\tilde\alpha; \hat\mF_i) =\inf\{t: \bP(v\leq t) \geq \tilde\alpha, v\sim   \hat\mF_i\}$. The empirical distribution $\hat\mF$ is discrete with mass $p_{ij}^H$ on $V_i$,  we can have an explicit expression for $Q(\tilde\alpha; \hat\mF_i)$:
\[
Q(\tilde\alpha; \hat\mF_i) = \left\{\begin{array}{ll}
\overline{V}_0,& \tilde\alpha =0,\\
\overline{V}_i,&  \sum_{j=1}^{i-1} p^H_{ij}< \tilde \alpha\leq  \sum_{j=1}^{i} p^H_{ij}, i = 1,\ldots, n, \\
\overline{V}_{n+1},&  \sum_{j=1}^{n} p^H_{ij}< \tilde \alpha.\\
\end{array}\right.
\]
Hence, $Q(\tilde\alpha; \hat\mF_i)$ is non-decreasing and right-continuous piece-wise constant on $\tilde\alpha$, and $v_{i}^*$ can only change its value at $\sum_{j=1}^{k} p^H_{ij}$ for $k=1,\ldots, n$. The same is true for $Q(\tilde\alpha; \hat\mF)$. 

\item[(ii)] Given $\tilde\alpha$, when increasing $V_{n+1}$ from $V_{n+1}=v'$  to $V_{n+1} = v$ for $v > v'$, the empirical distribution $\hat\mF_i(v)$ dominates the empirical distribution $\hat\mF_i(v')$ by construction: $\forall \tilde\alpha$, we have
\[
\bP(t\leq \tilde\alpha|t\sim \hat\mF_i(v'))\geq \bP(t\leq \tilde\alpha|t\sim \hat\mF_i(v)).
\]
As a result, $Q(\tilde\alpha; \hat\mF_i(v))$ is non-decreasing on $v$ for any given $\tilde\alpha$, for $i=1,\ldots, n+1$.
\item[(iii)]  Suppose that $v\in C_V(X_{n+1})$.  Let $\tilde\alpha\in \Gamma$ be the smallest value such that
\[
\sum^{n+1}_{i=1}\frac{1}{n+1}\mathbbm{1}_{V_i \leq  Q(\tilde\alpha; \hat\mF_i(v))}  \geq \alpha,
\]
by definition, we  have $v\leq Q(\tilde\alpha; \hat\mF)$. Now, we consider $V_{n+1} = v'$ for $v'\leq v$.  By the monotonicity of $ Q(\tilde\alpha; \hat\mF_i(v))$ on $\tilde\alpha$ and $v$ from Lemma \ref{lem:monotone} (i) and (ii), we must have $\tilde\alpha'\geq \tilde\alpha^*$ where $\tilde\alpha'\in \Gamma$ is the smallest value satisfying 
\[
\sum^{n+1}_{i=1}\frac{1}{n+1}\mathbbm{1}_{V_i \leq  Q(\tilde\alpha'; \hat\mF_i(v'))}  \geq \alpha,
\]
Hence,  we have $v' \leq v \leq Q(\tilde\alpha; \hat\mF) \leq Q(\tilde\alpha'; \hat \mF)$ and $v'$ is included in the PI.  This concludes our proof.
\end{itemize}
\end{proof}

\subsection{Proof of Lemma \ref{lem:conditional}}
\begin{proof}
Let $\sigma$ be a permutation of numbers $1,2,\ldots, n+1$.   We know that
\begin{align*}
P(\sigma_{n+1} = i|\mathcal{T}) &= \frac{\#\{\sigma:\sigma_{n+1} = i\}}{\sum^{n+1}_{j=1}\#\{\sigma:\sigma_{n+1} = j\}} = \frac{1}{n+1}.
\end{align*}
Set $\mX=\{X_1, \ldots ,X_{n+1}\}$ be the unordered set of the features.  Since the function $V(., \mZ) $ and the localizer $H(., ., \mX)$ are fixed functions conditional on $\mathcal{T}$, and $\tilde\alpha$ (can be random) is independent of the data conditional $\mathcal{T}$, we obtain
\begin{align}
\label{eq:lem2_eq1}
&\bP(V_{n+1} \leq Q(\tilde\alpha; \sum^{n+1}_{j=1} p^{H}_{n+1,j} \delta_{V_{j}}) |\mathcal{T}, \tilde\alpha) \notag\\
=& \sum^{n+1}_{i=1} P(\sigma_{n+1} = i|\mathcal{T})\mathbbm{1}_{\{V_{n+1}\leq v^*_{n+1}(\sigma)|\mathcal{T}, \sigma_{n+1} = i\}}\notag\\
=& \sum^{n+1}_{i=1} \frac{1}{n+1}\mathbbm{1}_{\{v_i\leq v^*_{n+1}(\sigma), \sigma_{n+1} = i\}}
\end{align}
where 
\[
v^*_{i}(\sigma) \coloneqq Q(\tilde\alpha; \sum^{n+1}_{j=1}p^{H}_{\sigma_i,\sigma_j}\delta_{v_{\sigma_j}})=Q(\tilde\alpha; \sum^{n+1}_{j=1} \frac{H( x_{\sigma_i}, x_{\sigma_j})}{\sum^{n+1}_{j'=1}H( x_{\sigma_i}, x_{\sigma_{j'}})}\delta_{v_{\sigma_{j}}} )
\] 
is the realization of $v^*_i\coloneqq Q(\tilde\alpha; \hat\mF_i)$ under permutation $\sigma$,  conditional on $\mathcal{T}$ and $\tilde\alpha$.  We immediately observe that,
\begin{equation}
\label{eq:lem2_eq2}
v^*_i(\sigma) = v^*_{\sigma_i}
\end{equation}
Combine (\ref{eq:lem2_eq1}) and (\ref{eq:lem2_eq2}), we obtain that $\bP\{V_{n+1} \leq v^*_{n+1} |\mathcal{T}, \tilde\alpha\} =  \sum^{n+1}_{i=1} \frac{1}{n+1}\mathbbm{1}_{\{v_i\leq Q(\tilde\alpha; \hat\mF_i)\}}$. Marginalize over $\tilde\alpha|\mathcal{T}$, we have
\[
\bP\{V_{n+1} \leq Q(\tilde\alpha; \hat\mF_{n+1}) |\mathcal{T}\}  = \bE\{ \sum^{n+1}_{i=1} \frac{1}{n+1}\mathbbm{1}_{\{v_i\leq Q(\tilde\alpha; \hat\mF_i)\}}|\mathcal{T}\}
\]
\end{proof}

\subsection{Proof of Lemma \ref{lem:alpha_bound}}
\begin{proof}
\;
\textbf{Part (i):}  We divide the space into non-overlapping subregions $A_k = \{x: (k-1)h_n \leq d(x_0, x ) < k h_n\}$.   Then, 
\[
B(x_0) = \sum_{i=1}^{n+1} \exp(-\frac{d(x_0, X_i)}{h_n})\geq \exp(-1)\sum_{i=1}^n \mathbbm{1}_{X_i\in A_1},
\]
and $\mathbbm{1}_{X_i\in A_1}$ follows a Bernoulli distribution with success probability $q_i \geq \frac{1}{L} nh_n^{\beta}$ according to Assumption \ref{ass:regularity1} (ii).  We can apply Chernoff Bounds to lower bound $B(x_0)$:
\[
\bP(\sum_{i=1}^n \mathbbm{1}_{X_i\in A_1}\leq \frac{nh_n^\beta}{2L})\leq \exp(-\frac{nh_n^\beta}{8L}),\qquad \bP(B(x_0)\leq \frac{nh_n^{\beta}}{2eL})\leq \exp(-\frac{n h_n^{\beta}}{8L}).
\]
Using the partitions  $\{A_j\}$ and Assumption \ref{ass:regularity1} (i):
\begin{align}
\label{eq:lem4_eq1}
\Delta(x_0)&\leq L\sum_{i=1}^n d(x_0,X_i)\exp(-\frac{d(x_0, X_i)}{h_n}),\notag\\
&\leq L h_n\exp(1)\sum_{k=1}^{\infty}k \sum_{i: X_i\in A_k} \exp(-k)\notag\\
& \leq \min_{k_0}\left\{L\exp(1) k_0 h_n\sum_{k\leq k_0} \sum_{i:X_i\in A_k}H(x_0, X_i)+Lh_n\exp(1)\sum_{k>k_0}\sum_{i:X_i\in A_k}k\exp(-k)\right\}\notag\\
& \leq  \min_{k_0}\left\{eLk_0 h_n B(x_0)+eLh_n k_0 \exp(-k_0) n\right\}\notag\\
& \leq  eL\beta\lceil \ln h_n^{-1}\rceil h_n \left(B(x_0)+n h_n^\beta\right),
\end{align}
where we have taken $k_0 = \beta\lceil \ln h_n^{-1}\rceil$ at the last step. Hence, there exists a constant $C>0$ such that
\[
\frac{\Delta(x_0)}{B(x_0)\vee (n h_n^{\beta})}\leq 2eL  \beta\lceil \ln h_n^{-1}\rceil h_n \leq C\ln (h_n^{-1})h_n,\; \mbox{for all }x_0\in[0,1]^p.
\]

\textbf{Part (ii):}  Set $Z_{ij} = \frac{H_{ij}}{B_i}\left(\mathbbm{1}_{V_j < V_i} - P_{V|X_j}(V_i)\right)$, and $R_i = \sum_{j\neq i}Z_{ij}$. By Hoeffding's lemma, the centered variable $Z_{ij}$ is sub-Gaussian with parameter $\nu_{ij} = \frac{H_{ij}}{2B_i}$ for all $i, j$ and $V_i$, e.g., for all $j\neq i$: 
\[
\bE[\exp(\lambda Z_{ij})|V_i] \leq \exp(\frac{\nu^2_{ij}\lambda}{2}),\;\mbox{ for all } \lambda\in \real.
\]
Hence, the weighted sum $R_{i} $ is sub-Gaussian with parameter $\nu_i = \sqrt{\sum_{j\neq i, j\leq n}\frac{H_{ij}^2}{4B_i^2}}\leq \sqrt{\frac{1}{4B_i}}$ (recall that $H_{ij}\leq 1$ and $B_i = \sum_{j=1}^{n+1} H_{ij}$). Combining it with the sub-Gaussian concentration results, we obtain that
\[
\bP(|R_i|\geq t|\mX, V_i)\leq  2\exp(-\frac{t^2}{2\nu_i^2})\leq   2\exp(-2t^2 B_i),\;\mbox{for all } V_i, i = 1,\ldots, n+1.
\]
Take $t =\sqrt{\frac{\ln n}{B_i}}$, we obtain the desired bound.
\end{proof}
\section{Choice of H}
\label{app:method_h}
\subsection{Estimation of the default distance}
Let $\mathcal{V}$ be the CV fold partitioning when learning $V$. We will estimate the spread by learning  $|V_i|$ for $V_i$ from the cross-validation step and $i=1,\ldots, n_0$:
\[
V_i\leftarrow \hat{V}^{-i}(X_i^0, Y_i^0),
\]
where $\hat V^{-i}$ is the score function learned using samples excluding $i$. 

The spread learning step is using the same CV partitioning $\mathcal{V}$. To learn the spread $\rho(X)$,  we consider minimizing the MSE with the response $\log(|V_i|+\overline{|V_i|})$, with $\overline{|V_i|}$ be the mean absolute value for $V_i$ across samples in $\mathcal{D}_0=\{Z_i^0=(X_i^0, Y_i^0), i=1,\ldots, n_0\}$.  This additional term $\overline{|V_i|}$ is added to reduce the influence of samples with very small empirical $|V_i|$. 

We do not claim that learning $\rho(X)$ in such a way is always a good choice. This is a reasonable choice the for regression score. However, for quantile regression score, $|V_i|$ is large around regions with both severe under-coverage  and over-coverage.  Despite this, we observe that the LCP ends up similarly as CP with a poorly chosen $\hat\rho(x)$ for the quantile regression score in our empirical studies.

Our estimated $\hat\rho$ is defined as $\hat\rho = \exp(\hat f(x))$ where $\hat f(x)$ is the estimated function from the learning step.  We let $\rho_i = \hat \rho^{-i}(X_i^0)$ be the estimated spread from the cross-validation step. Let $J\in \real^{n_0\times p}$ be the Jacobian  matrix with $J_{i,} = \frac{\partial \hat f^{-i}(X_i^0)}{\partial X_i^0}$.  Let $u_{\parallel}\in \real^{p\times p_0}$ and $u_{\perp}\in \real^{p\times (p-p_0)}$ be the top $p_0$ and the remaining right singular vectors, with $p_0$ be a small constant. By default, $p_0 = 1$. We form the projection matrix $P_{\parallel}$ and $P_{\perp}$ with $u_{\parallel}$ and  $u_{\perp}$:
\[
P_{\parallel} = u_{\parallel}u_{\parallel}^\top, \; P_{\perp} = u_{\perp}u_{\perp}^\top.
\] 
The final dissimilarity measure $d(x_1, x_2)$ is a weighted sum of the three components, and
\[
d(x_1, x_2) =  \frac{d_1(x_1,x_2)}{\sigma_2}+\frac{\left(\omega d_2(x_1, x_2)+(1-\omega)d_3(x_1, x_2)\right)}{\sigma_1},
\]
where $d_2(x_1,x_2)$, $d_3(x_1, x_2)$ are projected distances onto $P_{\parallel}$ and  $P_{\perp}$, and $d_1(x_1, x_2)$ are distance in the space of the learned spreading function $\hat\rho(x_1), \hat\rho(x_2)$ as described in Section \ref{subsec:method_h}: 
\begin{itemize}
\item  $d_1(x_1, x_2) = \|\hat\rho(x_1)-\hat\rho(x_2)\|_2$.
\item  $d_2(x_1, x_2) = \|P_{\parallel}(x_1 - x_2)\|_2$.
\item  $d_3(x_1, x_2) = \|P_{\perp}(x_1 - x_2)\|_2$.
\end{itemize}
We set $\omega$ and $\sigma_1$, $\sigma_2$ as following:
\begin{itemize}
\item  Let $\mu_{\parallel}\slash \mu_{\perp}$ be the mean of $d_2(X^0_i,X^0_j)$ or  $d_3(X^0_i,X^0_j)$ for $i\neq j$, then we let $w = \frac{\mu_{\perp}}{\mu_{\perp}+\mu_{\parallel}}$.
\item  We let $\sigma_1$ be the mean of $\left(\omega d_2(X^0_i,X^0_j)+(1-\omega)d_3(X^0_i,X^0_j)\right)$ and $\sigma_2$ be that mean of $d_1(X^0_i,X^0_j)$, using all pairs $i\neq j$ from $\mathcal{D}_0$.
\end{itemize}
\subsection{Empirical estimate of the objective}
We want to minimize a penalized  average length of finite PIs:
\begin{align*}
J(h) & = \mbox{Average PI}^{finite}\mbox{ length}+\lambda\times \mbox{Average conditional  PI}^{finite}\mbox{ length variability}\\
s.t. \;& \bP(\mbox{Infinite PI}) \leq \delta.
\end{align*}
Let $\bE_{X}f(X)$ denote the expectation of some function $f(.)$ with over $X$.  In this tuning section, we consider two specific types of $V(.)$: the scaled regression score and the scaled quantile score, and
\[
V(X, Y) = \frac{1}{\sigma(X)}|Y - f(X)|,
\]
or 
\[
V(X, Y) = \frac{1}{\sigma(X)}\max\{q_{lo}(X) - Y, Y-q_{hi}(X)\}.
\]
These two score classes will include the four scores considered in our numerical experiments. Let $k^*$ be the selected index from Lemma \ref{lem:alg1}, for this two classes of scores, the PI of $Y_{n+1}$ is constructed as
\[
C(X_{n+1}) = [f(X_{n+1})-\sigma(X_{n+1})\overline{V}_{k^*}, f(X_{n+1})+\sigma(X_{n+1})\overline{V}_{k^*}],
\]
or 
\[
C(X_{n+1}) = [q_{lo}(X_{n+1})-\sigma(X_{n+1})\overline{V}_{k^*}, q_{hi}(X_{n+1})+\sigma(X_{n+1})\overline{V}_{k^*}].
\]
In both cases, the length over the constructed PI of $Y_{n+1}$ is additive on $\sigma(X_{n+1})\overline{V}_{k^*}$, and hence, minimizing PI of $Y_{n+1}$ is equivalent to minimizing $\sigma(X_{n+1})\overline{V}_{k^*}$, and the conditional variability of the PI  is the same as the variability of $\sigma(X_{n+1})\overline{V}_{k^*}$ conditional on $X_{n+1}$. Hence, after omitting components that do not depend on $h$, we can express the terms in the above objective as 
\begin{itemize}
\item $\mbox{Average PI}^{finite}\mbox{ length}$: $ \bE_{Z_{1:n}, X_{n+1}}\left[\sigma(X_{n+1})\overline{V}_{k^*}|k^* \leq n\right]$. It depends on $Z_{1:n}, X_{n+1}$ as well as the tuning parameter $h$. (Recall that when $k^* = n+1$, $\overline{V}_{n+1} = \infty$. )
\item $\mbox{Average conditional  PI}^{finite}\mbox{ length variability}$:  $\sqrt{\bE_{Z_{1:n}, X_{n+1}}\left[\sigma(X_{n+1})^2\left(\overline{V}_{k^*}-\mu(X_{n+1})\right)^2|k^* \leq n\right]}$, where $\mu(X_{n+1}) = \bE_{Z_{1:n}}\left[\overline{V}_{k^*}|k^* \leq n, X_{n+1}\right]$ is the average length finite PI at $X_{n+1}$, marginalized over $Z_{1:n}$.
\item  Average percent of infinite PI:  $ \bP(k^*= n+1)$.
\end{itemize}
We estimate the above quantities with empirical estimates  using $\mathcal{D}_0$. As in the previous section, we consider the case where the function form $V(X, Y)$ is estimated by CV and $V_i\leftarrow \hat{V}^{-i}(X_i^0, Y_i^0)$.  For example, we want to construct the score function $V(X, Y) = |Y-f(X)|$ where $f(X)$ is the mean prediction function. Then, $ \hat V^{-i}(X_i^0, Y_i^0)$ is calculated as 
\[
 \hat V^{-i}(X_i^0, Y_i^0) = |Y^0_i-\hat f^{-k}(x)|,
\]
where $\hat f^{-k}(.)$ is the learned mean function using data excluding fold $k$ that includes sample $i$.  We also estimate the spreads and define the distance on $\mathcal{D}_0$ using the CV estimates.

Given the dissimilarity measure $d_{ij}$ for any pair $(X_i^0, X_j^0)$, and thus $H_{ij}=\exp(-\frac{d_{ij}}{h})$ for a given $h$, we estimate the empirical loss for $h\in \{h_1,\ldots, h_m\}$ as below:
\begin{itemize}
\item Estimation of average length and infinite PI probability:
\begin{itemize}
\item We subsample $\tilde n = (n+1)\wedge n_0$ samples without replacement from $\mathcal{D}_0$, let the set be $\mathcal{S}$ and construct PI for each  sample $i\in \mathcal{S}$ with a calibration set $\mathcal{S}\setminus\{i\}$. Let $L_i$ be the scaled length for the constructed PI (scaled by $\sigma(X_i^0)$).
\item The probability of having infinite PI is estimated as $C_1(h)=\frac{\#\{i\in \mathcal{S},L_i = \infty\}}{\tilde n}$, and the average finite PI length is estimated as $C_2(h)=\frac{\sum_{i\in \mathcal{S},L_i < \infty} L_i}{\#\{b:L_{ib}<\infty\}\vee 1}$.
\end{itemize}
The above estimates can be repeated for multiple times when $n_0$ is much larger than $(n+1)$.
\item Estimation of conditional variability:
\begin{itemize}
\item Repeat $B$ times the PI construction: for $b=1,\ldots, B$, we subsample $n $ samples with replacement from $\mathcal{D}_0$, and let the length of scaled PI of $V$ at $Z_i^0$ be $L_{ib}$ for $i=1,\ldots, n_0$.
\item Calculate the finite conditional mean as $\mu_i = \frac{\sum_{b:L_{ib}<\infty}L_{ib}}{\#\{b:L_{ib}<\infty\}\vee 1}$.
\item Calculate the conditional variance as $s_i =  \frac{\sum_{b:L_{ib}<\infty}(L_{ib}-\mu_i)^2}{\#\{b:L_{ib}<\infty\}\vee 1}$. 
\item The average conditional variability for PI with finite length is estimated as $C_{3,h} = \sqrt{\frac{\sum_i \left(\#\{b: L_{ib}<\infty\}\times s_i\right)}{\#\{(i,b): L_{ib}<\infty\}}}$
\end{itemize}
\end{itemize}
We take $h$ from the candidate set to minimize the empirical objective:
\[
h = \arg\min_{C_1(h)\leq \delta} \left(C_{2}(h)+\lambda C_3(h)\right).
\]

\end{document}